**Kolmogorov equations for evaluating the boundary hitting of degenerate diffusion with unsteady drift**

Hidekazu Yoshioka [a, *]

[a] Japan Advanced Institute of Science and Technology, 1-1 Asahidai, Nomi, Ishikawa 923-1292, Japan
[*] Corresponding author: yoshih@jaist.ac.jp, ORCID: 0000-0002-5293-3246

**Abstract**

Jacobi diffusion is a representative diffusion process whose solution is bounded in a domain under certain drift and diffusion coefficient conditions. However, the process without such conditions has not been thoroughly investigated. We explore a Jacobi diffusion whose drift coefficient is affected by another deterministic process, causing the process to hit the boundary of a domain in finite time. The Kolmogorov equation (a degenerate elliptic partial differential equation) for evaluating the boundary hitting of the proposed Jacobi diffusion is then presented and analyzed, with several conditional arguments, some of which are addressed computationally. We also investigate a related mean-field-type (McKean–Vlasov) self-consistent model arising in tourism management, where the drift depends on the index for sensor boundary hitting, thereby confining the process to a domain with higher probability. We propose a finite difference method for the linear and nonlinear Kolmogorov equations, which yields a unique numerical solution because of discrete ellipticity if the discount is positive. The accuracy of the finite difference method critically depends on the regularity of the boundary condition, and the use of high-order discretization is not always effective. Finally, we computationally investigate the mean field effect.

**CRediT Author Contributions:** All parts of this study were conducted by the author.

**Funding statement:** This study was supported by the Japan Society for the Promotion of Science (KAKENHI No. 22K14441 and 22H02456), the Japan Science and Technology Agency (PRESTRO, No. JPMJPR24KE), and Nippon Life Insurance Foundation (Environmental Research Grant for Young Researchers in 2024, No. 24).

**Data availability statement:** Data will be made available upon reasonable request to the corresponding author.

**Competing interests:** The author has no competing interests.

**Declaration of generative AI in scientific writing:** The author did not use generative AI technology to obtain the results in this manuscript.

**Permission to reproduce material from other sources:** N.A.

## 1. Introduction

### 1.1 Study background

Stochastic processes confined to domains arise in many problems in science and engineering. Lee and Whitmore [1] modeled disease progression as the hitting of a stochastic process to an upper boundary of a domain. Statistical evaluation of the first hitting time of stochastic processes is important for efficient molecular simulation [2]. Stability analysis of stochastic mechanical and structural systems has been performed in phase spaces [3,4]. Analyzing the boundary behavior of a stochastic process has thus been of great interest because of its relevance to a variety of applied studies. In this paper, we focus on diffusion processes, which are continuous-time stochastic processes driven by Brownian motion.

The diffusion process hits the boundary of a domain if the drift of the process is outward and/or if there is diffusion near the boundary, whereas a boundary cannot be touched by the process if the drift is inward and there is no diffusion (i.e., degenerate diffusion) near the boundary. Statistical evaluation of a stochastic process in a domain can be addressed by solving the associated Kolmogorov equation. Kolmogorov equations are second-order degenerate elliptic or parabolic partial differential equations (PDEs) whose solutions are conditional statistics of diffusion processes (e.g., Chapter 9 of Øksendal [5]). The boundary behavior of second-order degenerate elliptic PDEs has been systematically studied by Oleinik and Radkevič [6]; a mathematical tool, the Fichera function, is used to determine whether a boundary condition is necessary at each point of the boundary of a domain, depending on the drift and diffusion coefficients. Option pricing models have been computationally analyzed by effectively accounting for degenerate diffusion [7-9]. The hitting to a boundary has been analyzed in mechanical physics [10], fluid dynamics [11], chemistry [12], and population genetics as well [13].

A diffusion process can be confined to a domain by suitably choosing the coefficients of its governing stochastic differential equation (SDE). Systems of stochastic processes preserving certain ordering properties have been formulated by appropriately specifying the degeneration of diffusion coefficients [14] and jumps [15]. Sine–Wiener noise as a randomized sinusoidal function has been employed in the analysis of a potential well model [16] and a tumor growth model [17]. Among the existing bounded diffusion processes, the Jacobi diffusion process (sometimes called the Wright–Fisher process) (Chapter 6 of Alfonsi [18]) is the simplest. Its governing SDE is given by

$$\underbrace{\mathrm{d}X_t}_{\text{Increment}} = \underbrace{(b - aX_t)\mathrm{d}t}_{\text{Drift}} + \underbrace{c\sqrt{X_t(1-X_t)}\mathrm{d}B_t}_{\text{Diffusion}},\ t > 0 \tag{1}$$

subject to an initial condition $X_0 \in D \equiv (0,1)$. Here, $t \geq 0$ is time, $a,b,c > 0$ are constants, and $B = (B_t)_{t\geq 0}$ is a 1-D standard Brownian motion. The SDE (1) is understood in the Itô's sense. The parameters $a$, $b$, and $c$ represent the mean reversion, source, and noise intensity, respectively. The degenerate diffusion coefficient $\sqrt{X_t(1-X_t)}$ effectively bounds the solutions to the SDE (1) in $D$. More specifically, the process $X = (X_t)_{t\geq 0}$ is pathwise (i.e., in a strong sense) unique and is bounded in $D$ if and only if the following condition is satisfied (Proposition 6.2.1 of Alfonsi [18]):

$$2b \geq c^2 \text{ and } 2(a-b) \geq c^2 . \tag{2}$$

This condition means that the reversion $a$ is larger than the source $b$ and that the noise intensity $c$ is small. This kind of necessary and sufficient condition has been employed in a variety of problems, including but not limited to energy management [19,20], wind speed modeling [21], carbon emission markets [22], pairs trading [23], and random matrix theory in physics [24]. D'Onofrio [25] discussed the first hitting time of a Jacobi diffusion to a threshold placed interior to a domain.

If condition (2) is not satisfied, then the process $X$ hits the boundary of $D$ in a finite time with a probability of 1, making the statistical evaluation of the escape probability of the process from the domain and the analysis of boundary behavior meaningful. In population genetics, boundary mutation has been investigated, where the hitting of the Jacobi process to boundary points of a domain is of primary interest [26]. As a recent application of Jacobi processes, sustainable tourism management has been discussed, where the solutions represent travel demand, and the boundaries correspond to the states of over- and no-tourism [27,28]. In these studies, overtourism—characterized by pollution and the destruction of the environment and heritage due to the arrival of an excessive number of tourists at a tourism site [29-32]—has been regarded as a saturation of travel demand. However, to the best of the authors' knowledge, such boundary-hitting problems have been much less studied.

### 1.2 Aim and contribution

The aim of this paper is to formulate and analyze Jacobi diffusions driven by unsteady drifts, where their solutions hit a boundary of the domain in finite time. More specifically, we assume that the source $b$ is time dependent and that the second condition (2) is not initially satisfied but becomes satisfied later. Thus, the boundary stability of the process changes over time. Such a problem has recently been investigated numerically in Yoshioka [28] on the basis solely of a Monte Carlo simulation without deep mathematical analysis, but we address it through Kolmogorov equations as a novel approach. Unsteady Jacobi diffusion arises in the tourism management problem, where the hitting of the upper boundary of a domain represents the state of overtourism [27,28]. The proposed system of SDEs conceptually describes the transient dynamics of tourism demand subject to an external force to mitigate the state of overtourism: the hitting to the boundary of a solution.

The use of Kolmogorov equations is more advantageous than Monte Carlo simulations because the former evaluates the boundary-hitting phenomenon for any initial condition at once [e.g., 33-35], whereas the latter needs to be run for each initial condition. In contrast, a disadvantage of the approach based on Kolmogorov equations is the curse of dimensionality, as computing a high-dimensional problem requires exceptionally large computational resources such as memory and time; however, this does not apply to our case because it is a two-dimensional problem. Moreover, Monte Carlo simulations become inefficient if the target process is of a mean-field type, where the drift or diffusion coefficient of the target SDE depends on its law. We argue that Kolmogorov equations can be formulated in a mean-field case as well, focusing on a Jacobi process of a mean-field type.

In this paper, we address two Kolmogorov equations. The first Kolmogorov equation of our Jacobi diffusion with unsteady drift is a degenerate elliptic PDE whose boundary condition should be prescribed along a part of the boundary (Chapter 1 of Oleinik and Radkevič [6]). More specifically, the boundary condition of the equation is prescribed where the underlying Jacobi process probably occurs (**Figure 1** in **Section 2**). Solutions to this Kolmogorov equation should be understood in a viscosity sense under a certain continuity assumption along the boundary [36,37], and we explicitly provide a necessary condition for this continuity. We present some sort of mathematical analysis results about the first Kolmogorov equation, including a comparison result based on a conditional argument about subsolution and supersolution. In particular, the uniqueness is conditional on an additional comparison-type ordering property, and that the numerical section (**Section 4**) is best understood as approximating the specific value function suggested by the stochastic representation rather than a fully characterized unique viscosity solution in the usual PDE sense.

The second Kolmogorov equation investigated in this paper has a nonlinear drift coefficient, representing mean-field-type (McKean–Vlasov) self-consistent feedback [38]. The second Kolmogorov equation is reduced to the first equation when nonlinearity is omitted; hence, the first equation is generalized. The second Kolmogorov equation includes nonlinearity such that its solution appears in the drift coefficient. This kind of nonlinearity is usually not encountered in classical stochastic control [5] but rather in the master equations of mean field games [38,39]. In the context of tourism modeling, this situation implies that a central planner, such as a government, adaptively controls tourism dynamics to avoid overtourism in a feedback manner by considering the statistical evaluation of the hitting of the SDE to a boundary. We investigate the second Kolmogorov equation computationally; therefore, this paper does not provide a well-posedness theory, comparison principle, or convergence analysis for the equation. Nevertheless, we believe that this equation motivates a new strategy for modeling tourism.

Our model is related to reach-avoid problems where a Kolmogorov equation or its controlled version—a Bellman equation—governs the hitting or related probability of a system to a preferred or unpreferred state [40-42]. A key difference between their problems and ours is the diffusion coefficient and boundary conditions; the former addresses deterministic or jump-driven cases, strictly elliptic problems with nonvanishing diffusion, or situations where the entire domain is reachable by the system. In contrast, our problem targets systems driven by degenerate diffusion, where only a part of the boundary is reachable.

For both Kolmogorov equations, the regularity of boundary conditions is crucial in analyzing their well-posedness. This point is numerically investigated on the basis of a finite difference method of the monotone type [43]. A monotone finite difference method is often convergent in the sense of viscosity solutions and is theoretically desirable, but its computational performance is low, with a convergence rate of one or lower [e.g., 44]. We therefore additionally examine a higher-order nonmonotone finite difference method of the filtered scheme that adaptively blends low- and high-order discretization depending on the regularity of the numerical solutions [45]. This scheme has been successfully applied to degenerate parabolic PDEs [46,47] and nonlinear first-order PDEs [48,49]; however, its application to the problem of our form—a multidimensional degenerate elliptic problem with partial boundary conditions—has not been addressed. Our discretization method is nonlinear and (degenerate) elliptic [e.g., 50,51] and yields a unique numerical solution for both linear and nonlinear Kolmogorov equations if a positive discount exists.

The absence of analytical solutions to the Kolmogorov equations poses an obstacle to validating the finite difference method, which can be resolved at least for the first Kolmogorov equation by exploiting its stochastic representation to which a Monte Carlo simulation applies. By using a fine Monte Carlo result as the ground truth, we examine the convergence rate of the finite difference method with or without the filtered scheme. Applying the filtered scheme proves to be effective because the convergence rate exceeds one when the boundary condition is continuously differentiable, suggesting that the boundary regularity critically affects the performance of the finite difference methods. In particular, a convergence order of less than one is observed for the case where the first Kolmogorov equation governs the hitting probability of the underlying Jacobi diffusion to the boundary to which a discontinuous boundary condition applies. The second Kolmogorov equation is also numerically computed to investigate the influence of the mean-field effect. By relying on computational techniques to complement mathematical analysis, this paper contributes providing a new perspective on degenerate elliptic PDEs related to (mean-field) Jacobi diffusion processes and their numerical analysis.

The rest of this paper is organized as follows. **Section 2** introduces classical Jacobi diffusion and its variant with unsteady drift. The mean-field-type self-consistent version is also introduced in this section. The Kolmogorov equations associated with these models are presented and analyzed as well. **Section 3** presents a finite difference method to be applied to the Kolmogorov equations. **Section 4** addresses numerical investigations of solutions to the Kolmogorov equations. **Section 5** summarizes this study and presents its perspectives. The **Appendix** contains proofs (**Section A.1**) and auxiliary results (**Section A.2**).

## 2. Jacobi diffusion process and Kolmogorov equation

### 2.1 Jacobi diffusion

We work on a filtered probability space $\left(\Omega,\mathbb{F},\left(\mathbb{F}_t\right)_{t\geq 0},\mathbb{P}\right)$ (e.g., [52]). The classical Jacobi diffusion $X=\left(X_t\right)_{t\geq 0}$ is a continuous-time stochastic process governed by the Itô's SDE:

$$\mathrm{d}X_t=\left(b-aX_t\right)\mathrm{d}t+c\sqrt{X_t\left(1-X_t\right)}\mathrm{d}B_t,\ 0<t<\tau \tag{3}$$

subject to an initial condition $X_0\in D=\left(0,1\right)$, where $\tau$ is the first hitting time of $X$ to the boundaries of $D$:

$$\tau=\inf\left\{t>0\middle|X_t=0\text{ or }1\right\}. \tag{4}$$

We have $\tau=+\infty$ under (2), and then the SDE (3) admits a unique, strong solution (Proposition 6.2.1 of Alfonsi [18]). In contrast, if (2) is not satisfied, then $X$ probably hits the boundaries of $D$ in a finite time.

In the sequel, without any loss of generality, we assume the nondimensionalization $a=1$.

### 2.2 Jacobi diffusion with unsteady drift

We are interested in Jacobi diffusion with unsteady drift. The condition (2) is crucial for determining the boundary stability of the classical Jacobi diffusion, and our interest is in the case where this condition is initially violated but becomes satisfied later. Such a problem is encountered in modeling tourism [27-28]. An attractive tourist destination can attract an excessive number of visitors, triggering overtourism and degrading or destroying the local environment and heritage. Examples include Mt. Fuji in Japan [30] and various world heritage sites [32]. In our context, $X$ represents tourism demand at a particular destination, and $Y$ represents the net utility, including the costs and benefits that travelers gain by arriving there. We consider that excessively large values of $Y$ trigger overtourism, which is modeled as the hitting of $X$ to the upper boundary of $D$. We also assume that some policy measures are implemented to reduce the values of $Y$, making the hitting event to the boundary less possible. We show that this dynamic can be described by Jacobi diffusion with unsteady drift.

Modeling this transition is addressed by replacing $b$ with some process $Y = (Y_t)_{t \geq 0}$. To investigate various timescales of the transition in a unified way, we employ a logistic-type equation following the literature on rate-induced tipping [e.g., 53,54]:

$$\mathrm{d}Y_t = RY_t\left(1 - \frac{Y_t}{\overline{Y}}\right)\mathrm{d}t,\ t > 0 \tag{5}$$

subject to an initial condition $Y_0$, where $R > 0$ is the transition rate and $\overline{Y} > 0$ is the asymptotic source rate ( $\lim_{t \to +\infty} Y_t = \overline{Y}$ ). We assume that $Y_0 > \overline{Y}$. Under this setting, the process gradually decreases from $Y_0$ to $\overline{Y}$ following a sigmoidal profile whose transition width is $O\left(R^{-1}\right)$. The transition from $Y_0$ to $\overline{Y}$ becomes faster as $R$ increases.

We consider the following system of SDEs where $X$ is a Jacobi diffusion with unsteady drift with $\tau$ being given by (4):

$$\mathrm{d}\begin{pmatrix} X_t \\ Y_t \end{pmatrix} = \begin{pmatrix} Y_t - X_t \\ RY_t\left(1 - \frac{Y_t}{\overline{Y}}\right) \end{pmatrix}\mathrm{d}t + \begin{pmatrix} c\sqrt{X_t\left(1 - X_t\right)} \\ 0 \end{pmatrix}\mathrm{d}B_t,\ 0 < t < \tau. \tag{6}$$

We focus on the hitting of $X$ to the upper boundary $\{x = 1\}$. Considering (2), we assume the following conditions so that the desired boundary stability is obtained:

$$2\overline{Y} > c^2,\ Y_0 > 1,\ 2\left(1 - \overline{Y}\right) > c^2. \tag{7}$$

We assume that $\overline{Y} < 1 < Y_0$. The first condition of (7) means that the hitting of $X$ to the lower boundary $\{x = 0\}$ of $D$ never occurs. The second and third conditions then mean that the hitting of $X$ to the upper boundary $\{x = 1\}$ of $D$ occurs near the initial time. More specifically, the hitting to the upper boundary never occurs after the time $t = t_c$ satisfying $Y_{t_c} = 1 - \frac{c^2}{2}$. Such a $t_c$ exists uniquely since $\overline{Y} < 1 - \frac{c^2}{2} < Y_0$.

Without significant loss of generality, we focus on the specific case where $\overline{Y} = 1 - \delta$ and $Y_0 = 1+\delta$ with $\delta \in (0,1)$. Afterward, we introduce the decreasing normalized process $Z = (Z_t)_{t\geq 0}$ with the initial condition $Z_0 = 1$ and the asymptotic value $\lim_{t\to+\infty} Z_t = 0$:

$$Z_t = \frac{Y_t - (1-\delta)}{1+\delta-(1-\delta)} = \frac{Y_t-(1-\delta)}{2\delta}. \tag{8}$$

Substituting (8) into (6) yields the following nondimensional system:

$$\mathrm{d}\begin{pmatrix} X_t \\ Z_t \end{pmatrix} = \begin{pmatrix} 2\delta Z_t + 1 - \delta - X_t \\ -2\delta R\left(\dfrac{2\delta}{1-\delta} Z_t + 1\right) Z_t \end{pmatrix} \mathrm{d}t + \begin{pmatrix} c\sqrt{X_t(1-X_t)} \\ 0 \end{pmatrix} \mathrm{d}B_t,\ 0 < t < \tau. \tag{9}$$

With respect to the unique existence of strong (pathwise) solutions to the system (9), we have the following proposition. Its proof uses Theorem 1 in Yamada and Watanabe [55].

***Proposition 1***

*Under the assumption (7), a unique strong solution $(X,Z)$ exists in the system (9) up to time $\tau$.*

By **Proposition 1**, in what follows, we always understand the process $(X,Z)$ after $\tau$ as $(X_t,Z_t) = (X_\tau, Z_\tau)$ ($t \geq \tau$). We also show that under the assumption (7), it holds true that

$$\tau = \inf\{t>0 | X_t = 0 \text{ or } 1\} = \inf\{t > 0 | X_t = 1\}, \tag{10}$$

excluding the possibility of hitting the process $X$ to the lower boundary of $D$. Its proof uses the comparison result of Theorem 1.1 in Ikeda and Watanabe [56].

***Proposition 2***

*Under the assumption (7), the relationship (10) holds true. Moreover, the system (9) admits a unique strong solution.*

We show that $\tau = +\infty$ if $Z_0 \leq \rho \equiv \frac{1}{2} - \frac{c^2}{4\delta}$. Namely, boundary hitting never occurs if the source rate is small.

***Proposition 3***

*It follows that $\tau = +\infty$ if $Z_0 \leq \rho$.*

***Remark 1*** One can study the hitting of $X$ to the lower boundary of $D$ by the following transformation $X_t \to 1 - X_t$ along with suitably redefining the coefficients.

### 2.3 Kolmogorov equation

Boundary hitting of the process $X$ governed by (9) can be evaluated through the conditional expectation

$$V(x,z)=\mathbb{E}^{x,z}\left[f(Z_\tau)e^{-\eta\tau}\mathbb{I}(\tau<+\infty)\right],\ (x,z)\in D\times D. \tag{11}$$

Here, $\mathbb{I}(\tau<+\infty)$ is the indicator function for the event $\{\tau<+\infty\}$ ( $\mathbb{I}(\tau<+\infty)=1$ if $\tau<+\infty$ and $\mathbb{I}(\tau<+\infty)=0$ otherwise), $\mathbb{E}^{x,z}$ is the expectation conditioned on $(X_0,Z_0)=(x,z)$, $f:\bar{D}\to\mathbb{R}$ is a bounded and Borel measurable function, $\bar{D}=[0,1]$ is the closure of $D$, and $\eta\geq 0$ is the discount rate with which the expectation $V$ focuses more on the events near future as $\eta$ increases. Our target is how much probability sample paths of $X$ hits the upper boundary of $D$, namely, the escape probability from $D$. In this case, we may take $\eta=0$ and $f\equiv 1$. Another example is a smooth $f$ such that $f(\rho)=0$ (meaning that this condition becomes clearer later).

It is important to observe that $\mathbb{I}(\tau<+\infty)$ in (11) is actually redundant if $\eta>0$ since

$$\begin{aligned}\mathbb{E}^{x,z}\left[f(Z_\tau)e^{-\eta\tau}\right]&=\mathbb{E}^{x,z}\left[f(Z_\tau)e^{-\eta\tau}\mathbb{I}(\tau=+\infty)\right]+\mathbb{E}^{x,z}\left[f(Z_\tau)e^{-\eta\tau}\mathbb{I}(\tau<+\infty)\right]\\&=\mathbb{E}^{x,z}\left[f(Z_\tau)e^{-\eta\tau}\mathbb{I}(\tau<+\infty)\right]\end{aligned},\ (x,z)\in D\times D. \tag{12}$$

The formula for the stochastic representation of the conditional expectation of the form (11) (e.g., Theorem 15.3.1 of Pascucci [57]) tells us that $V$ if it is sufficiently smooth, solves the Kolmogorov equation

$$\eta V=(2\delta z+1-\delta-x)\frac{\partial V}{\partial x}-2\delta R\left(\frac{2\delta}{1-\delta}z+1\right)z\frac{\partial V}{\partial z}+\frac{1}{2}c^2x(1-x)\frac{\partial^2 V}{\partial x^2},\ (x,z)\in D\times D \tag{13}$$

subject to the boundary condition

$$V=f(z),\ (x,z)\in\Gamma, \tag{14}$$

where $\Gamma=\{x=1\}\times(\rho,1)$ is a part of the boundary to which the process $X$ probably hits (**Figure 1**).

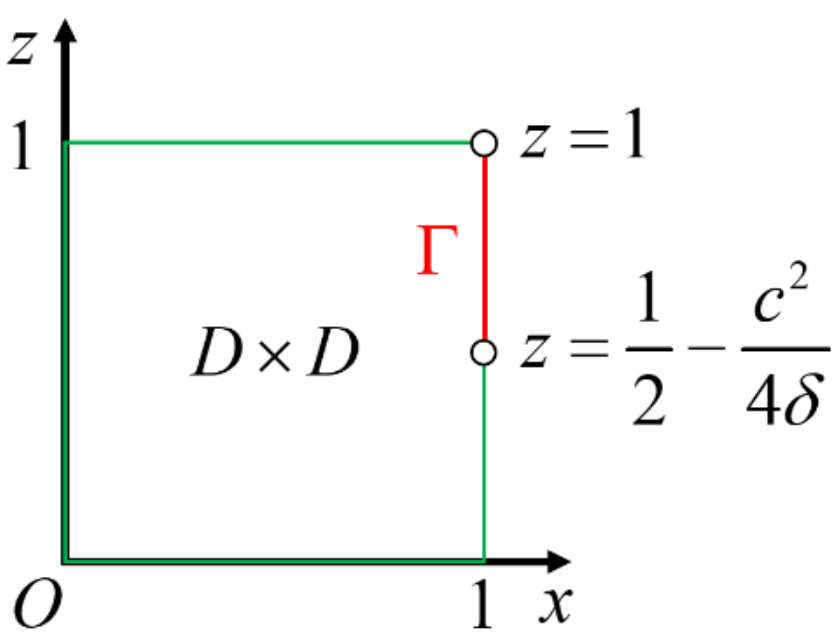

**Figure 1.** The domain of the Kolmogorov equation. The boundary condition is prescribed only along $\Gamma$.

As shown in (14), the boundary condition is prescribed only on $\Gamma$, corresponding to the fact that boundary hitting occurs only if the value of the process $Z$ is large (**Proposition 3**). Intuitively, this is because the characteristic curves corresponding to the drift of the system (9) are inward, the diffusion degenerates at $x=0,1$, and there is no diffusion in the $z$ direction. Mathematically, the partial boundary condition is justified by Fichera theory, which states that a Dirichlet boundary condition should be prescribed where the Fichera function—computable given the boundary geometry and the drift and diffusion coefficients along the boundary (e.g., Eq. (1.1.3) of Oleinik and Radkevič [6])—is negative or where the diffusion perpendicular to the boundary is positive.

Let $\mathbf{n}=(n_1,n_2)$ be the inward normal vector along the boundaries: $\mathbf{n}=(1,0)$ on $x=0$, $\mathbf{n}=(-1,0)$ on $x=1$, $\mathbf{n}=(0,1)$ on $z=0$, and $\mathbf{n}=(0,-1)$ on $z=1$. The drift $\mathbf{b}=(b_1,b_2)$ and (symmetric) diffusion coefficients $\mathbf{d}=\left[d_{i,j}\right]_{1\le i,j\le 2}$ of the PDE (13) are given as follows:

$$\mathbf{b}=\left(2\delta z+1-\delta-x,-2\delta R\left(\frac{2\delta}{1-\delta}z+1\right)z\right) \text{ and } \mathbf{d}=\begin{pmatrix}\frac{1}{2}c^2x(1-x) & 0\\ 0 & 0\end{pmatrix}. \tag{15}$$

Then, the Fichera function $\mathrm{Fic}(x,z)$ is given as follows (we write $x_1=x$ and $x_2=z$ here):

$$\mathrm{Fic}(x,z)=(n_1,n_2)\cdot\left(b_1-\frac{\partial}{\partial x_1}D_{11}-\frac{\partial}{\partial x_1}D_{12},b_2-\frac{\partial}{\partial x_1}D_{21}-\frac{\partial}{\partial x_2}D_{22}\right), \tag{16}$$

which in our case is given by

$$\mathrm{Fic}(x,z)=\left((2\delta z+1-\delta-x)-\frac{1}{2}c^2(1-2x)\right)n_1-2\delta R\left(\frac{2\delta}{1-\delta}z+1\right)zn_2. \tag{17}$$

Boundary values should be specified when $\mathrm{Fic}(x,z)<0$. Now, elementary calculations give

$$\mathrm{Fic}(x,z)=\begin{cases}2\delta z+1-\delta-\frac{1}{2}c^2 & (\text{along } x=0)\\ -2\delta z+\delta-\frac{1}{2}c^2 & (\text{along } x=1)\\ 0 & (\text{along } z=0)\\ 2\delta R\frac{1+\delta}{1-\delta} & (\text{along } z=1)\end{cases}, \tag{18}$$

yielding (14). More specifically, the first line in (18) for $x=0$ is positive due to the first condition in (7) that yields $\overline{Y}=1-\delta>c^2/2$. The second line for $x=1$ is negative for $z>\rho$, stating that the boundary value should prescribed when for $x=1$ and $z>\rho$. The third and fourth lines for $z=0,1$ are trivially not negative.

A difficulty in dealing with the Kolmogorov equation (13) can be illustrated by the following example. Assume that $f(z)=1$ and $\eta=0$ with which $V$ of (11). Then, there are no escape events when $Z_0<\rho$; hence, we must have $V(\cdot,z)=0$ for $z<\rho$. In contrast, the Kolmogorov equation (13) admits a constant function $V\equiv 1$ as a classical solution that complies with the boundary condition, which is clearly different from the escape probability. This implies that when $\eta=0$, there exists more than one "solutions" to the Kolmogorov equation (13). Assuming that $\eta>0$ is crucial in the context of the variational approach for degenerate elliptic PDEs in weighted Sobolev spaces: Chapters 1.5 and 1.6 of Oleinik and Radkevič [6] for generic degenerate elliptic PDEs and Feehan and Pop [58] for a specific case in finance.

In summary, the boundary condition should be prescribed only along $\Gamma = \{x = 1\} \times (\rho, 1)$, and a reasonable value of $V$ along $\{x = 1\} \times (0, \rho]$ is 0. The boundary value of $V$ is therefore continuous only if $f(\rho) = 0$. Similarly, the boundary value of $V$ is continuously differentiable only if $f(\rho) = 0$ and $\left.\frac{\mathrm{d}f(z)}{\mathrm{d}z}\right|_{z \to \rho - 0} = 0$. In this view, specifying $f = 1$ in $(\rho, 1)$, which is a natural choice for evaluating the hitting probability of $X$, is an irregular case that is challenging in numerical computation.

As explained above, understanding the meaning of "solution" to the Kolmogorov equation (13) is important because it is a degenerate elliptic PDE whose solutions are generally neither continuous nor differentiable [36]. In particular, boundary conditions are necessary only along $\Gamma$. A suitable solution for such cases is the viscosity solution for generalized boundary value problems: a weak solution that does not necessarily satisfy the boundary condition along the entire boundary. In the sequel, $C(\Theta)$, $LSC(\Theta)$, $USC(\Theta)$, and $C^{\infty}(\Theta)$ represent the collections of functions that are continuous, lower-semicontinuous, upper-semicontinuous, arbitrary times continuously differentiable in a domain $\Theta$, respectively.

**Definition 1** is a notion of viscosity solutions suited to our Kolmogorov equation (e.g., Definition 1 of Rokhlin [37]). For any generic sufficiently regular $\varphi : \bar{D} \times \bar{D} \to \mathbb{R}$ and any $(x, z) \in \bar{D} \times \bar{D}$, we set

$$\begin{aligned} \mathbb{L}\varphi(x,z) = {} & \eta\varphi(x,z) - (2\delta z + 1 - \delta - x)\frac{\partial \varphi(x,z)}{\partial x} + 2\delta R\left(\frac{2\delta}{1-\delta}z + 1\right) z \frac{\partial \varphi(x,z)}{\partial z} \\ & - \frac{1}{2}c^2 x(1-x)\frac{\partial^2 \varphi(x,z)}{\partial x^2} \end{aligned} . \quad (19)$$

***Definition 1***

*A function $\phi \in USC(\bar{D} \times \bar{D})$ is said to be a viscosity subsolution if for any $\varphi \in C^{\infty}(\bar{D} \times \bar{D})$ such that $\phi - \varphi$ on $\bar{D} \times \bar{D}$ is locally strictly maximized at $(x_0, z_0)$ then*

$$\mathbb{L}\varphi(x_0, z_0) \le 0,\ (x_0, z_0) \in D \times D, \quad (20)$$

$$\min\{\mathbb{L}\varphi(x_0, z_0), \varphi(x_0, z_0) - f(z_0)\} \le 0,\ (x_0, z_0) \in \partial(D \times D). \quad (21)$$

*A function $\phi \in LSC(\bar{D} \times \bar{D})$ is said to be a viscosity supersolution if for any $\varphi \in C^{\infty}(\bar{D} \times \bar{D})$ such that $\phi - \varphi$ on $\bar{D} \times \bar{D}$ is locally strictly minimized at $(x_0, z_0)$ then*

$$\mathbb{L}\varphi(x_0, z_0) \ge 0,\ (x_0, z_0) \in D \times D, \quad (22)$$

$$\max\{\mathbb{L}\varphi(x_0, z_0), \varphi(x_0, z_0) - f(z_0)\} \ge 0,\ (x_0, z_0) \in \partial(D \times D). \quad (23)$$

*A function $\phi \in C(\bar{D} \times \bar{D})$ is said to be a viscosity solution if it is a viscosity supersolution and is a viscosity subsolution.*

To clarify the mathematical mechanism of **Definition 1** and reconcile it with the Fichera analysis, we explicitly state how the partial boundary-value problem is realized. Although the boundary operator involving $f(z_0)$ is formally written over the entire boundary $\partial(D \times D)$ in equations (21) and (23), the Dirichlet data $f$ is effectively imposed only on $\Gamma$. This is a fundamental feature of the relaxed generalized boundary conditions for degenerate equations.

On the segment $\Gamma = \{x = 1\} \times (\rho, 1)$, the Fichera function is strictly negative ($\mathrm{Fic} < 0$) due to the outward drift, which prevents sub- and supersolutions from satisfying the internal operator $\mathbb{L}\phi = 0$ up to the boundary. Thus, the $\min/\max$ conditions strictly enforce $\phi = f(z_0)$ on $\Gamma$.

Conversely, on the complementary boundary $\partial(D \times D) \setminus \Gamma$, the Fichera function is non-negative ($\mathrm{Fic} \geq 0$). This segment acts as a state-constraint (or outflow) boundary, meaning that sample paths of the underlying diffusion process (9) cannot hit or cross it. Mechanistically in (21) and (23), because the characteristics point inward or the diffusion perpendicular to the boundary vanishes, the internal PDE inequalities $\mathbb{L}\varphi \leq 0$ (for subsolutions) and $\mathbb{L}\varphi \geq 0$ (for supersolutions) are naturally satisfied by the test functions at the boundary. Consequently, the $f$-involved terms are automatically bypassed by the $\min/\max$ selection, ensuring that no boundary conditions are active on $\partial(D \times D) \setminus \Gamma$.

Finally, the value 0 on the portion $\{x = 1, z \leq \rho\}$ is not an independently prescribed boundary condition, but a natural deterministic and probabilistic consequence of **Proposition 3**; since the process never reaches this portion, the value 0 is naturally selected as the unique bounded solution compatible with the stochastic representation (11).

In view of **Definition 1**, if $f \in C(\bar{D})$, then $V$ is a viscosity solution to the Kolmogorov equation (13), as shown below. The proof is an application of the stochastic Perron's method for generalized boundary value problems [37] combined with the pathwise uniqueness (**Proposition 2**) of the process $(X, Z)$ up to the first hitting time $\tau$. This methodology has been applied to a finite-fuel control problem [59]. Now, **Proposition 4** addresses the existence of a viscosity solution to the Kolmogorov equation (13).

***Proposition 4***

*If $f \in C(\bar{D})$ and $\eta > 0$, then there exists a viscosity subsolution $\underline{V}$ and a viscosity supersolution $\bar{V}$ to the Kolmogorov equation (13). Moreover, it follows that $V$ is a viscosity solution to the Kolmogorov equation (13) and that $V$ in (11) is bounded as $\underline{V} \leq V \leq \bar{V}$ on $\bar{D} \times \bar{D}$.*

**Proposition 4** does not address the uniqueness of viscosity solutions. Under an additional condition, $V$ becomes the unique viscosity solution to the Kolmogorov equation (13). One may hope that this is due to a continuity assumption of viscosity solutions along the boundary of the domain $D$ with the comparison theorem: Theorem 5.1 of Koike [60].

***Proposition 5***

*Assume that* $f \in C(\bar{D})$ *and* $\eta > 0$. *If for any viscosity subsolution* $\underline{V}'$ *and viscosity supersolution* $\bar{V}'$, *it follows that* $\underline{V}' \le \bar{V}'$ *on* $\bar{D} \times \bar{D}$*; then,* $V$ *of (11) is the unique viscosity solution to (13).*

We should comment on the continuity assumption ($\underline{V}' \le \bar{V}'$ on $\bar{D} \times \bar{D}$) in **Proposition 5**. Initially, we considered that this condition would follow from a direct application of the comparison theorem for generalized boundary value problems (Theorem 5.1 of Koike, [60]). However, it may fail because of the lack of suitable regularity of the diffusion coefficient, as discussed in **Section A.2 of the Appendix** (if there is no diffusion, this issue does not arise, but it is a less interesting case). By applying the argument based on the Fichera function, we can analytically derive a nontrivial necessary condition for the boundary continuity of viscosity solutions: $f \in C(\bar{D})$ should satisfy $f(\rho) = 0$ and $f(z) = 0$ for $0 \le z \le \rho$. The criticality of this assumption is investigated numerically in **Section 4**. Finally, the assumptions made in **Proposition 5**, particularly those on subsolution and supersolutions, are currently conditional.

***Remark 2*** One may consider that a generalized Itô's formula (e.g., Chapter 11.3.2 of Bensoussan [61]) combined with the continuity of the process $(X, Z)$ up to $\tau$ may apply to $V$ under a suitable condition of boundary data (e.g., Proof of Theorem 1.3 of Feehan and Pop [58]). However, this seems to be difficult for our case because [58] assumes that the diffusion term does not degenerate inside the domain, whereas our case does; it has no diffusion in the $z$ direction. This implies that the regularity of $V$, such as $V \in C_{\mathrm{loc}}\left((D \times D) \cup \Gamma\right) \cap C^{1,1}(D \times D)$ (i.e., $V$ is locally continuous in $(D \times D) \cup \Gamma$ and is continuously differentiable with Lipschitz continuous derivatives in $D \times D$), would not hold true.

### 2.4 Mean-field-type (McKean–Vlasov) self-consistent model

The mean-field-type self-consistent model in this paper is a system of SDEs whose coefficient depends on $V$ [e.g., 38,39,62]. We focus on the case where the dynamics of $Z$ and hence that of $X$ depend on $V$:

$$\mathrm{d}\begin{pmatrix} X_t \\ Z_t \end{pmatrix} = \begin{pmatrix} 2\delta Z_t + 1 - \delta - X_t \\ -2\delta R\left(\dfrac{2\delta}{1-\delta} Z_t + 1\right) Z_t \omega\left(V\left(X_t, Z_t\right)\right) \end{pmatrix} \mathrm{d}t + \begin{pmatrix} c\sqrt{X_t\left(1 - X_t\right)} \\ 0 \end{pmatrix} \mathrm{d}B_t,\ 0 < t < \tau \tag{24}$$

with a bounded continuous function $\omega : \mathbb{R} \to (0, +\infty)$. The SDE (24) represents a situation where the system dynamics are dynamically regulated by $V$, which is understood here as an objective to be reduced, suggesting a feedback mechanism from the objective to the system. In view of (18), this modification does not change the boundary condition of the Kolmogorov equation. The system (24) implicitly depends on its law through $V$, which is a conditional expectation.

The Kolmogorov equation in the mean-field case is the following nonlinear degenerate elliptic

PDE subject to the boundary condition (14):

$$\eta V=\left(2\delta z+1-\delta-x\right)\frac{\partial V}{\partial x}-2\delta R\left(\frac{2\delta}{1-\delta}z+1\right)z\omega\left(V\right)\frac{\partial V}{\partial z}+\frac{1}{2}c^{2}x\left(1-x\right)\frac{\partial^{2}V}{\partial x^{2}},\ \left(x,z\right)\in D\times D. \quad (25)$$

For simplicity, we focus on the following specific $\omega$ unless otherwise specified:

$$\omega\left(v\right)=\frac{1}{2}\left(1+\tanh\left(\frac{v}{\kappa}\right)\right),\ v\in\mathbb{R} \quad (26)$$

with a constant $\kappa>0$. This choice corresponds to the promotion of the decay of $Z_t$ as $V\left(X_t,Z_t\right)$ increases, representing feedback between objective $V$ and dynamics. We can return to the linear case by setting $\omega\equiv 1$.

Computationally, we can also use different $\omega$ as demonstrated in **Section 4**, but we mainly investigate (26) because it is smooth and bounded, contributing to the unique existence of numerical solutions to the finite difference method (see the proof of **Proposition 4** in **Section 3**). One may also consider a negative $\omega$, which represents the situation where the increase of $Z_t$ is suppressed as $V\left(X_t,Z_t\right)$ increases. However, in the context of tourism modeling, this situation is less interesting because it potentially leads to an earlier hitting to the boundary and, consequently, an earlier realization of overtourism with a higher probability.

Finally, one may consider viscosity solutions to the nonlinear Kolmogorov equation for the linear case. This is theoretically possible, but its uniqueness and existence would be more difficult issues, which are worth investigating in future works.

## 3. Numerical method

### 3.1 Computational grid

We use a computational grid consisting of uniformly placed vertices $\mathrm{P}_{i,j}$ ($i,j=0,1,2,...$) whose locations are $\left(x_i,z_j\right)=\left(i/N,j/N\right)$, with $N\in\mathbb{N}$ being the computational resolution. The domain $D\times D$ is thus uniformly discretized by these vertices with the interval $h=1/N$. The numerical approximation of $V$ at $\mathrm{P}_{i,j}$ is denoted by $V_{i,j}$. The goal of our finite difference method is to obtain $V_{i,j}$ at all the vertices. The set of vertices $\Gamma_{\mathrm{num}}$ corresponding to $\Gamma$ is given by $\Gamma_{\mathrm{num}}=\left\{\mathrm{P}_{N,j}\middle|z_j>\rho,\ j=0,1,2,...,N\right\}$. We set $V_{-1,\cdot}=V_{N+1,\cdot}=V_{\cdot,-1}=0$. Our numerical method addresses the cases $\eta>0$ and $\eta=0$ in a unified way, but the treatment of the latter needs case because the latter case admits multiple solutions; see the next section.

### 3.2 Monotone finite difference method

We first present a monotone finite difference method. At each vertex exterior to $\Gamma_{\mathrm{num}}$, each term of the Kolmogorov equation (25) is discretized as follows (for the linear version (13), we simply set $\omega\equiv 1$):

$$\eta V\to\eta V_{i,j}, \quad (27)$$

$$\frac{1}{2}c^2 x(1-x)\frac{\partial^2 V}{\partial x^2} \to \frac{1}{2}c^2 x_i(1-x_i)\frac{V_{i+1,j}-2V_{i,j}-V_{i-1,j}}{h^2}, \tag{28}$$

$$(2\delta z+1-\delta-x)\frac{\partial V}{\partial x} \to (2\delta z_j+1-\delta-x_i)\begin{cases}\dfrac{V_{i+1,j}-V_{i,j}}{h} & (2\delta z_j+1-\delta-x_i \geq 0)\\[2ex] \dfrac{V_{i,j}-V_{i,j-1}}{h} & (2\delta z_j+1-\delta-x_i < 0)\end{cases}, \tag{29}$$

$$-2\delta R\left(\frac{2\delta}{1-\delta}z+1\right)z\omega(V)\frac{\partial V}{\partial z} \to -2\delta R\left(\frac{2\delta}{1-\delta}z_j+1\right)z_j\bar{\omega}(V_{i,j},V_{i,j-1})\frac{V_{i,j}-V_{i,j-1}}{h}. \tag{30}$$

Here, $\bar{\omega}:\mathbb{R}\times\mathbb{R}\to(0,+\infty)$ is a function with $\bar{\omega}(v,v)=\omega(v)$ specified later. At each vertex in $\Gamma_{\mathrm{num}}$, we prescribe the boundary condition

$$V_{i,j}=f(z_j). \tag{31}$$

Collecting all (27)-(31) yields the system

$$-A_{i,j}V_{i-1,j}+(\eta+A_{i,j}+B_{i,j}+C_{i,j}+D_{i,j})V_{i,j}-B_{i,j}V_{i+1,j}=C_{i,j}V_{i,j-1}+E_{i,j} \tag{32}$$

for $i,j=0,1,2,...,N$. The coefficients $A_{i,j},B_{i,j},C_{i,j},D_{i,j},E_{i,j}$ are given as follows: at each vertex exterior to $\Gamma_{\mathrm{num}}$,

$$A_{i,j}=\frac{1}{2}c^2 x_i(1-x_i)\frac{1}{h^2}-(2\delta z_j+1-\delta-x_i)\frac{1}{h}\times\begin{cases}0 & (2\delta z_j+1-\delta-x_i\geq 0)\\ 1 & (2\delta z_j+1-\delta-x_i<0)\end{cases}, \tag{33}$$

$$B_{i,j}=\frac{1}{2}c^2 x_i(1-x_i)\frac{1}{h^2}+(2\delta z_j+1-\delta-x_i)\frac{1}{h}\times\begin{cases}1 & (2\delta z_j+1-\delta-x_i\geq 0)\\ 0 & (2\delta z_j+1-\delta-x_i<0)\end{cases}, \tag{34}$$

$$C_{i,j}=2\delta R\left(\frac{2\delta}{1-\delta}z_j+1\right)z_j\bar{\omega}(V_{i,j},V_{i,j-1})\frac{1}{h}, \tag{35}$$

and $D_{i,j}=E_{i,j}=0$, and at each vertex in $\Gamma_{\mathrm{num}}$, we have

$$D_{i,j}=1,\ E_{i,j}=f(z_j),\ A_{i,j}=B_{i,j}=C_{i,j}=0. \tag{36}$$

We should design $\bar{\omega}$ so that the discretized system becomes degenerate elliptic and hence satisfies the discrete maximum principle (e.g., Section 3.4 of Bonnans et al. [50] and Section 2.2 of Oberman [51]). We propose

$$\bar{\omega}(V_{i,j},V_{i,j-1})=\frac{1}{2}\left(1+\tanh\left(\frac{V_{i,j}}{\kappa}\right)\right)\mathbb{I}(V_{i,j}\geq V_{i,j-1})+\frac{1}{2}\left(1+\tanh\left(\frac{V_{i,j-1}}{\kappa}\right)\right)\mathbb{I}(V_{i,j}<V_{i,j-1}). \tag{37}$$

For any $v\in\mathbb{R}$, we have

$$\bar{\omega}(v,v)=\frac{1}{2}\left(1+\tanh\left(\frac{v}{\kappa}\right)\right)\mathbb{I}(v\geq v)+\frac{1}{2}\left(1+\tanh\left(\frac{v}{\kappa}\right)\right)\mathbb{I}(v<v)=\omega(v). \tag{38}$$

We check the degenerate ellipticity of the discretization (32) with (37). To see this, we rewrite (32) as follows by separating the linear and nonlinear parts:

$$\Xi_{i,j} \equiv \underbrace{-A_{i,j}V_{i-1,j} + \left(\eta + A_{i,j} + B_{i,j} + D_{i,j}\right)V_{i,j} - B_{i,j}V_{i+1,j} - E_{i,j}}_{\text{Linear part}} + \underbrace{\Lambda_{i,j}}_{\text{Nonlinear part}} = 0 , \tag{39}$$

where

$$\begin{aligned}\Lambda_{i,j} &= 2\delta R\left(\frac{2\delta}{1-\delta}z_j + 1\right)z_j\bar{\omega}\left(V_{i,j},V_{i,j-1}\right)\frac{V_{i,j}-V_{i,j-1}}{h} \\ &= 2\delta R\left(\frac{2\delta}{1-\delta}z_j + 1\right)z_j \\ &\times\left\{\frac{1}{2}\left(1+\tanh\left(\frac{V_{i,j}}{\kappa}\right)\right)\max\left\{\frac{V_{i,j}-V_{i,j-1}}{h},0\right\}+\frac{1}{2}\left(1+\tanh\left(\frac{V_{i,j-1}}{\kappa}\right)\right)\min\left\{\frac{V_{i,j}-V_{i,j-1}}{h},0\right\}\right\}\end{aligned} . \tag{40}$$

We have the following proposition showing that the system (32) admits a unique numerical solution $\left\{V_{i,j}\right\}_{i,j=0,1,2,\ldots,N}$ because of the degenerate ellipticity (Section 2.2 of Oberman [51]). Here, the discretized equation is said to be degenerate elliptic if $\Xi_{i,j}$ is nondecreasing with respect to $V_{i,j}$ and is nonincreasing with respect to all $\left\{V_{k,l}\right\}_{(k,l)\neq(i,j)}$ for all $i,j=0,1,2,\ldots,N$. The discretized equation is said to be elliptic if "nondecreasing" in the previous sentence is strengthened to "strictly increasing."

***Proposition 6***

*Assume that $\eta > 0$. Then, the system (32) is elliptic and admits a unique numerical solution $\left\{V_{i,j}\right\}_{i,j=0,1,2,\ldots,N}$.*

The assumption of positive discount $\eta > 0$ is essential in **Proposition 6** so that our discretized system becomes elliptic. Without that, the system is only degenerate elliptic.

Computationally, the particular form of the system (32) allows us to find a numerical solution from $j=0$ to $j=N$ in a cascading manner, and at each $i$, we only need to solve the tridiagonal system. We use a relaxation method from Oberman and Salvador [45], as explained in **Section 4**. The sweep direction is from $i=0$ to $i=N$.

For the linear case with $\omega \equiv 1$, owing to **Proposition 6**, if $\eta > 0$, then the numerical solution is $V_{i,j}=0$ for $i=0,1,2,\ldots$ and $j=0,1,2,\ldots\bar{j}-1$, where $\bar{j}$ is the largest integer $j<N$ such that $\mathrm{P}_{N,j}$ is exterior to $\Gamma_{\mathrm{num}}$. This follows from an induction argument starting from $j=0$ to a larger $j>0$ along with the ellipticity of the discretization scheme owing to $\eta>0$. Therefore, if $f$ is uniformly continuous for $\rho \leq z \leq 1$ such that $f\left(\rho\right)=0$, then we expect better convergence behavior of the numerical solutions. In contrast, if $\eta=0$, then for any real constant $\beta$, $V_{i,j}=\beta$ for $i=0,1,2,\ldots$ and $j=0,1,2,\ldots\bar{j}-1$ solves the system; hence, the numerical solutions are not uniquely determined. However, case $\beta=0$ is the most natural case because $\tau=+\infty$ when $Z_0<\rho$, as shown in **Proposition 3**. To resolve this nonuniqueness, we use the initial guess such that any numerical solution equals zero at all vertices, with which the solution with $\beta=0$ is obtained accordingly.

***Remark 3*** For the nonlinear case with $\omega \neq 1$, the existence of viscosity solutions to the Kolmogorov equation (25) is nontrivial. This is due mainly to the unusual form of the Hamiltonian associated with this equation, in our context, particularly the product term $\omega(V)\frac{\partial V}{\partial x}$ that resembles Burgers-type equations [38,39,63].

***Remark 4*** Alasseur et al. [39] applied upwind discretization to the nonlinear term of their Kolmogorov equation, which adaptively uses the information from each side depending on the sign of the nonlinear drift coefficient. We do not use their method to avoid the loss of the cascading nature of the discretized system, which can be solved from $j = 0$ to $j = N$. The filtered scheme introduced in the next subsection does not lose this advantage.

***Remark 5*** The unique results of the numerical solutions presented in **Proposition 6** also apply to the linear case ( $\omega = 1$ ) if $\eta > 0$.

### 3.3 Filtered scheme

We now present the finite difference method equipped with the filtered scheme. Intuitively, this is a nonmonotone finite difference method derived by adding a correction term to the monotone term presented in the earlier subsection. The added term is expected to improve the convergence rate, although there may be a loss of computational stability. For later use, we set the filter function

$$F(y) = \begin{cases} y & (-1 \leq y \leq 1) \\ 0 & (y > 1,\ y < -1) \end{cases}, \tag{41}$$

which is a discontinuous function that equals the identity mapping for small $|x|$ and has been used in the existing filtered schemes [e.g., 46,48].

The finite difference method equipped with the filtered scheme is designed as follows at each vertex exterior to $\Gamma_{\text{num}}$ (the discretization at each vertex in $\Gamma_{\text{num}}$ remains unchanged):

$$-A_{i,j}V_{i-1,j} + \left(\eta + A_{i,j} + B_{i,j} + C_{i,j} + D_{i,j}\right)V_{i,j} - B_{i,j}V_{i+1,j} = C_{i,j}V_{i,j-1} + E_{i,j} + G_{i,j}, \tag{42}$$

where the coefficient $G_{i,j}$ corresponds to the term added by the filtered scheme:

$$G_{i,j} = G_{i,j,x} + G_{i,j,z}, \tag{43}$$

where

$$G_{i,j,x} = \begin{cases} F\left(\dfrac{1}{h^{1/2}}\left(2\delta z_j + 1 - \delta - x_i\right)\left(\Delta V_{i,j,1,x} - \Delta V_{i,j,2,x}\right)\right) & (3 \leq i \leq N-3) \\ 0 & (\text{Otherwise}) \end{cases} \tag{44}$$

and

$$G_{i,j,x} = \begin{cases} F\left(\dfrac{1}{h^{1/2}}\left(-2\delta R\left(\dfrac{2\delta}{1-\delta}z_j + 1\right)z_j\bar{\omega}\left(V_{i,j}, V_{i,j-1}\right)\right)\left(\Delta V_{i,j,1,z} - \Delta V_{i,j,2,z}\right)\right) & (j \geq 3) \\ 0 & (\text{Otherwise}) \end{cases}. \tag{45}$$

Here, we set the first-order one-sided differences

$$\Delta V_{i,j,1,x} = \begin{cases} \dfrac{V_{i+1,j} - V_{i,j}}{h} & \left(2\delta z_j + 1 - \delta - x_i \geq 0\right) \\ \dfrac{V_{i,j} - V_{i,j-1}}{h} & \left(2\delta z_j + 1 - \delta - x_i < 0\right) \end{cases} \tag{46}$$

and

$$\Delta V_{i,j,1,z} = \frac{V_{i,j} - V_{i,j-1}}{h} . \tag{47}$$

We also set third-order one-sided differences (Section 2.3.8 of Festa et al. [64]).

$$\Delta V_{i,j,2,x} = \begin{cases} \dfrac{1}{h}\left(-\dfrac{11}{6}V_{i,j} + 3V_{i+1,j} - \dfrac{3}{2}V_{i+2,j} + \dfrac{1}{3}V_{i+3,j}\right) & \left(2\delta z_j + 1 - \delta - x_i \geq 0\right) \\ \dfrac{1}{h}\left(\dfrac{11}{6}V_{i,j} - 3V_{i-1,j} + \dfrac{3}{2}V_{i-2,j} - \dfrac{1}{3}V_{i-3,j}\right) & \left(2\delta z_j + 1 - \delta - x_i < 0\right) \end{cases} \tag{48}$$

and

$$\Delta V_{i,j,2,z} = \frac{1}{h}\left(\frac{11}{6}V_{i,j} - 3V_{i,j-1} + \frac{3}{2}V_{i,j-2} - \frac{1}{3}V_{i,j-3}\right). \tag{49}$$

The term $G_{i,j}$ is the difference between the first- and third-order differences and serves as the correction of the discretization, which is activated if their difference is sufficiently small. The scaling factor $h^{1/2}$ comes from a formal convergence argument between low- and high-order discretization methods $\left|\Delta V_{i,j,1,x} - \Delta V_{i,j,2,x}\right|, \left|\Delta V_{i,j,1,z} - \Delta V_{i,j,2,z}\right| = O\left(h^{1/2}\right)$ (Remark 1 of Oberman and Salvador [45]), and we follow it. The filtered scheme is not applied to the discretization near boundaries.

According to Theorem 1 of Oberman and Salvador [45], numerical solutions generated by the proposed finite difference method are stable (i.e., bounded irrespective of the resolution $N$). Moreover, our computational results in **Section 4** agree with those of the Monte Carlo simulations. The last assumption is nontrivial because of the existence of the coefficient $G_{i,j}$, which is nonlinear with respect to the numerical solution. As shown in **Section 4**, our numerical solutions are bounded irrespective of computational resolution.

***Remark*** *7* In view of Theorem 1 of Oberman and Salvador [45], the direct use of a high-order discretization (setting $F(y) = y$) should be avoided. In contrast, the monotone scheme (setting $F \equiv 0$) is within the application range of this theorem. We therefore do not examine the direct use of a high-order discretization.

## 4. Numerical computation

### 4.1 Monte Carlo simulation

We numerically analyze the Kolmogorov equations for which analytical solutions have not been found. Therefore, we validate the finite difference method against the “ground truth” computed by the naïve Monte

Carlo simulation based on the Euler–Maruyama method. The total number of sample paths is $2\times10^{6}$, and the time increment for temporal integration is $5\times10^{-6}$, with the total number of time steps $2\times10^{6}$.

We can numerically estimate the quantity $V$ by computing sample paths multiple times, although this can be computationally intensive. A drawback of the Monte Carlo approach is that sample paths must be iterated for each initial condition $(x,z)$. In contrast, a Kolmogorov equation provides a solution at all points in the domain $D\times D$. Another drawback of the Monte Carlo method is its extreme inefficiency for computing $V$ when the drift depends on $V$, as in (25), because calculating the system itself requires information from $V$; concurrently finding this information poses a challenging task.

### 4.2 Boundary regularity

We first analyze the linear Kolmogorov equation by comparing the finite difference method with and without the filtered scheme against the Monte Carlo results. Here, we set the following parameter values that satisfy (7): $\delta=0.5$, $c=0.4$, and $R=0.2$. This parameter set is referred to as the nominal case in the sequel.

**Tables 1-2** compare the finite difference and Monte Carlo results with and without the filtered scheme where $V(1,\cdot)$ is discontinuous at $z=\rho$, where we set $f(z)=f_1(z)=1$ so that $V$ represents the escape probability of $X$ from $D$. Similarly, **Tables 3-4** compare the results where $V(1,\cdot)$ is Lipschitz continuous but not continuously differentiable at $z=\rho$: $f(z)=f_2(z)=\max\left\{1,10\left(2\delta z-\delta+c^2/2\right)\right\}$. **Tables 5-6** compare the results where $V(1,\cdot)$ is continuously differentiable at $z=\rho$: $f(z)=f_3(z)=\max\left\{2\delta z-\delta+c^2/2,0\right\}^2$. In all the cases above, we examine both with discount $\eta=0.1$ and without it $\eta=0$. We use the following relaxation method to improve the convergence (for the monotone scheme, we simply omit $G_{i,j}$):

$$V_{i,j}\Rightarrow wV_{i,j}+(1-w)\frac{A_{i,j}V_{i-1,j}+B_{i,j}V_{i+1,j}+C_{i,j}V_{i,j-1}+E_{i,j}+G_{i,j}}{\eta+A_{i,j}+B_{i,j}+C_{i,j}+D_{i,j}} \tag{50}$$

with the relaxation factor $w\in[0,1)$ where "$\Rightarrow$" indicates the point update. Here, a larger $w$ yields a more stable iteration but with a slower convergence. We set $w=1/2$. The sweeping at each $j$ is judged to be converged when the absolute difference between the new and old numerical solution values become smaller than $10^{-12}$ at all $i=0,1,2,...,N$. To reduce computational costs, we check the convergence at each one hundred iterations so that the actual error threshold is more severe.

**Figures 2** shows the computed $V$ for each $f$ with $\eta=0$; those with $\eta>0$ have similar profiles but with smaller values and are therefore not presented. As shown in **Figure 2**, the numerical solutions do not exhibit artificial oscillations. For $f=f_2$, **Figure 3** compares the Monte Carlo result, the numerical solution of the monotone finite difference method, and that with the filtered scheme along $z=1$, demonstrating that the numerical solution enhanced with the filtered scheme better captures the Monte Carlo result and is less diffusive. The convergence of the finite difference method is typically achieved within $O(10^2)$ to $O(10^3)$ iterations at each $j$ (see also **Section A3** in **Appendix**).

**Tables 1-6** imply that the convergence rate improves as the regularity of $f$ increases, where the convergence rate is smaller than 0.5 for $f_1$ in both the monotone and filtered schemes, and the error of the latter is slightly smaller than that of the former. The convergence rate is lower with the filtered scheme than without it in some cases, which is attributed to the discretization becoming less sensitive as the resolution increases. Therefore, the benefits of using the filtered scheme are particularly significant under lower computational resolution. According to **Tables 5-6**, the filtered scheme outperforms the monotone scheme in terms of both the error and the convergence rate for a sufficiently regular $f$.

For $f_2$ and $f_3$, the convergence rates of the monotone and filtered schemes improve and exceed 1 in the filtered scheme, particularly with the continuously differentiable $f = f_3$. Nevertheless, the convergence rate remains smaller than 2. One may expect that the convergence rate should approach 2 because the diffusion term is discretized by a centered difference and the drift terms by filtered schemes, which possibly has third-order accuracy; however, this did not occur in our case. We consider that this discrepancy is due to the low regularity of the current Kolmogorov equation, which does not include any diffusion in the $z$ direction and has low regularity near the boundary $x = 1$, particularly when $f = f_1$ is discontinuous. The other $f$ has higher regularities, but a local loss of regularity in the solutions may occur near the boundary, especially around $z = \rho$. The solution inside $D \times D$ (except for $z < \rho$) is influenced by the boundary data; therefore, both the loss of regularity and computational errors propagate from the boundary into the interior of the domain.

The computational results also suggest that the absence of a discount decreases the formal convergence rate of the proposed finite difference method when it is applied to a Kolmogorov equation, at least in the present case. In an extreme case (not addressed in this paper), a solution to a Kolmogorov equation for a problem without a discount can become unbounded. For example, if $\eta = 0$, then the mean residence time of the process $(X, Z)$ will diverge along $z = \rho$ and is $+\infty$ if $z < \rho$ because the process $X$ never hits the boundary when $Z_0$ is small. Therefore, designing a statistic $V$ that is bounded and hence computable is critical for analyzing the behavior of degenerate diffusion processes such as Jacobi diffusion.

We also conduct a sensitivity analysis with respect to parameter values that potentially affect the system dynamics, as shown in **Figure 4**, where we examine cases with a smaller $R = 0.05$, greater noise intensity $c = 0.6$, and greater variation $\delta = 0.75$. Here, we use $f_3$. The computational resolution is $N = 800$. The smaller $R$ leads to $V$ that varies less with respect to $x$ because of the greater possibility of the boundary hitting $X$. A greater noise intensity $c$ results in a greater probability of hitting the boundary at earlier times. The greater variation $\delta$ between the maximum and minimum $Y$ results in the greater likelihood of the process $X$ reaching the boundary when it is relatively large. From a tourism management standpoint, the computational results imply that a prolonged increase in tourist flow would occur because of slower policy-making or more unpredictable tourist flow.

**Table 1.** Comparison between the finite difference and Monte Carlo results: without filtered scheme: $f = f_1$. The convergence rate at resolution $N$ is calculated as $\log_2\left(\text{Error}_N / \text{Error}_{2N}\right)$, where $\text{Error}_N$ is the error at resolution $N$. The same applies to the other tables in this paper.

| $N$ | $\eta = 0$ (Ground truth: 0.27156) | | $\eta = 0.1$ (Ground truth: 0.23843) | |
|---|---|---|---|---|
| | Error | Convergence rate | Error | Convergence rate |
| 200 | 0.07207 | 0.517 | 0.06229 | 0.528 |
| 400 | 0.05037 | 0.501 | 0.04322 | 0.510 |
| 800 | 0.03559 | 0.478 | 0.03035 | 0.289 |
| 1600 | 0.02555 | 0.454 | 0.02484 | 0.655 |
| 3200 | 0.01866 | | 0.01578 | |

**Table 2.** Comparison between the finite difference and Monte Carlo results: with filtered scheme: $f = f_1$.

| $N$ | $\eta = 0$ (Ground truth: 0.27156) | | $\eta = 0.1$ (Ground truth: 0.23843) | |
|---|---|---|---|---|
| | Error | Convergence rate | Error | Convergence rate |
| 200 | 0.04862 | 0.383 | 0.04133 | 0.390 |
| 400 | 0.03729 | 0.387 | 0.03155 | 0.390 |
| 800 | 0.02852 | 0.386 | 0.02407 | 0.389 |
| 1600 | 0.02182 | 0.336 | 0.01838 | 0.388 |
| 3200 | 0.01729 | | 0.01405 | |

**Table 3.** Comparison between the finite difference and Monte Carlo results: without filtered scheme: $f = f_2$.

| $N$ | $\eta = 0$ (Ground truth: 0.22502) | | $\eta = 0.1$ (Ground truth: 0.19889) | |
|---|---|---|---|---|
| | Error | Convergence rate | Error | Convergence rate |
| 200 | 0.05121 | 0.760 | 0.04493 | 0.767 |
| 400 | 0.03024 | 0.751 | 0.02640 | 0.428 |
| 800 | 0.01796 | 0.724 | 0.01962 | 1.061 |
| 1600 | 0.01087 | 0.690 | 0.00941 | 0.697 |
| 3200 | 0.00674 | | 0.00580 | |

**Table 4.** Comparison between the finite difference and Monte Carlo results: with filtered scheme: $f = f_2$.

| $N$ | $\eta = 0$ (Ground truth: 0.22502) | | $\eta = 0.1$ (Ground truth: 0.19889) | |
|---|---|---|---|---|
| | Error | Convergence rate | Error | Convergence rate |
| 200 | 0.01920 | 0.610 | 0.01671 | 0.622 |
| 400 | 0.01258 | 0.530 | 0.01086 | 0.991 |
| 800 | 0.00871 | 0.512 | 0.00547 | 0.068 |
| 1600 | 0.00611 | 0.502 | 0.00522 | 0.508 |
| 3200 | 0.00431 | | 0.00367 | |

**Table 5.** Comparison between the finite difference and Monte Carlo results: without filtered scheme: $f = f_3$.

| $N$ | $\eta = 0$ (Ground truth: 0.42869) | | $\eta = 0.1$ (Ground truth: 0.38441) | |
|---|---|---|---|---|
| | Error | Convergence rate | Error | Convergence rate |
| 200 | 0.21658 | 1.000 | 0.19478 | 1.002 |
| 400 | 0.10833 | 0.984 | 0.09726 | 0.986 |
| 800 | 0.05477 | 0.969 | 0.04910 | 0.971 |
| 1600 | 0.02797 | 0.955 | 0.02505 | 0.957 |
| 3200 | 0.01443 | | 0.01291 | |

**Table 6.** Comparison between the finite difference and Monte Carlo results: with filtered scheme: $f = f_3$.

| $N$ | $\eta = 0$ (Ground truth: 0.42869) | | $\eta = 0.1$ (Ground truth: 0.38441) | |
|---|---|---|---|---|
| | Error | Convergence rate | Error | Convergence rate |
| 200 | 0.16443 | 1.335 | 0.14340 | 1.353 |
| 400 | 0.06520 | 1.429 | 0.05635 | 1.460 |
| 800 | 0.02422 | 1.564 | 0.02049 | 1.536 |
| 1600 | 0.00819 | 1.299 | 0.00707 | 1.245 |
| 3200 | 0.00333 | | 0.00298 | |

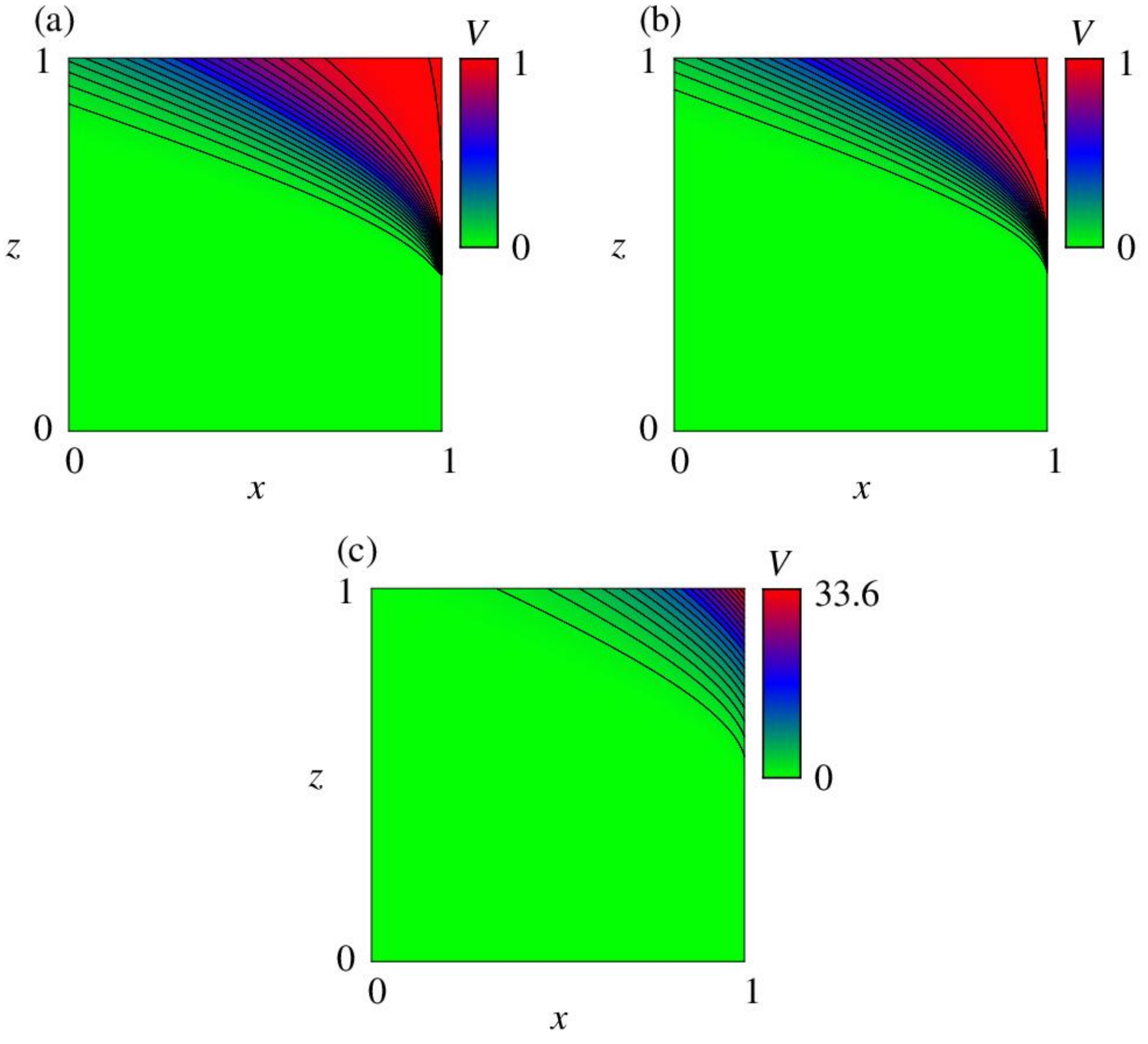

**Figure 2.** Computed $V$ with the filtered scheme for each $f$ with the filtered scheme: **(a)** $f_1$, **(b)** $f_2$, and **(c)** $f_3$. The computational resolution in the finite difference method is $N = 800$.

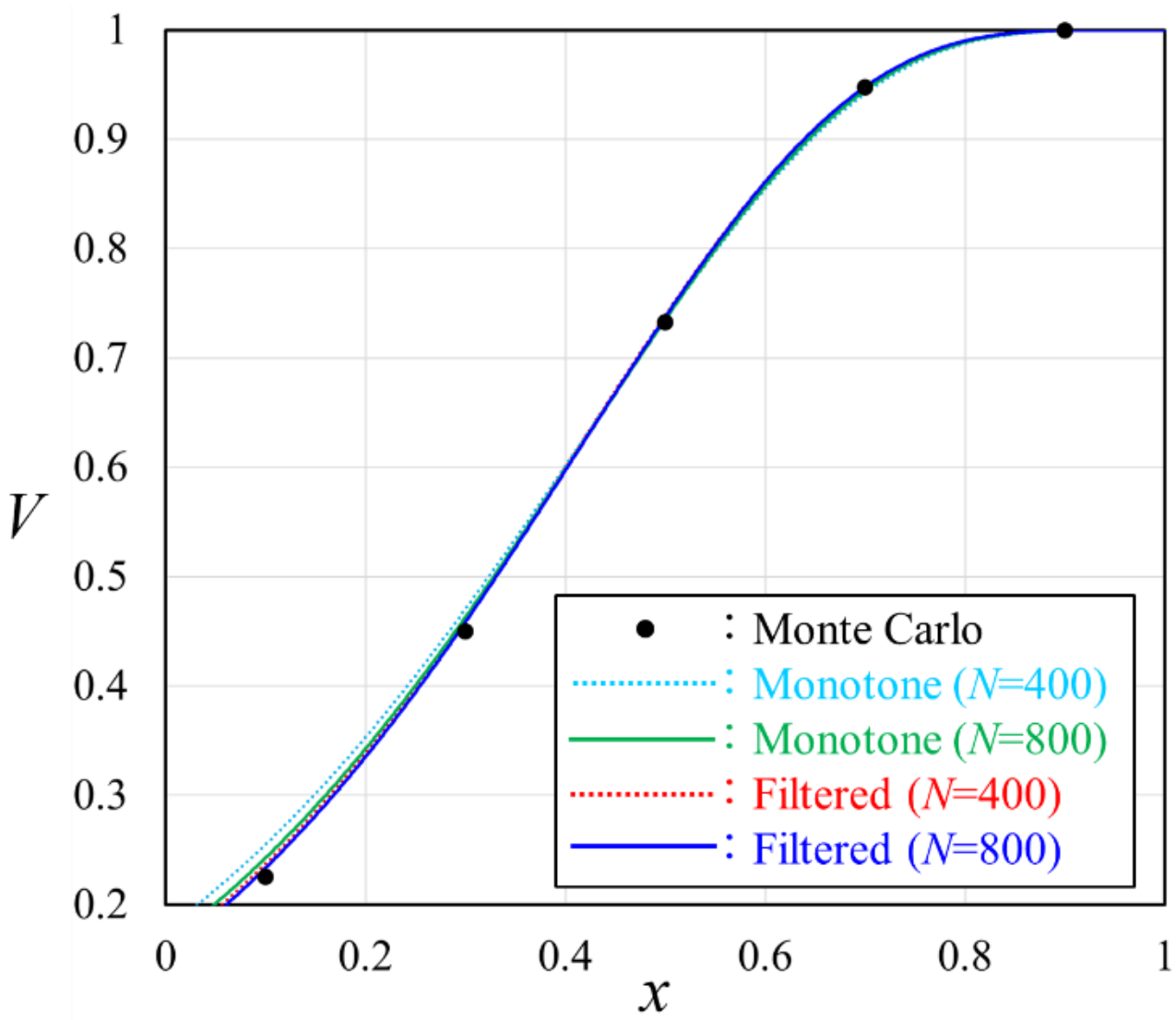

**Figure 3.** Monte Carlo result (circles), numerical solution of the monotone finite difference method (red curve), and that with the filtered scheme (blue curve) along $z = 1$ for $f = f_2$.

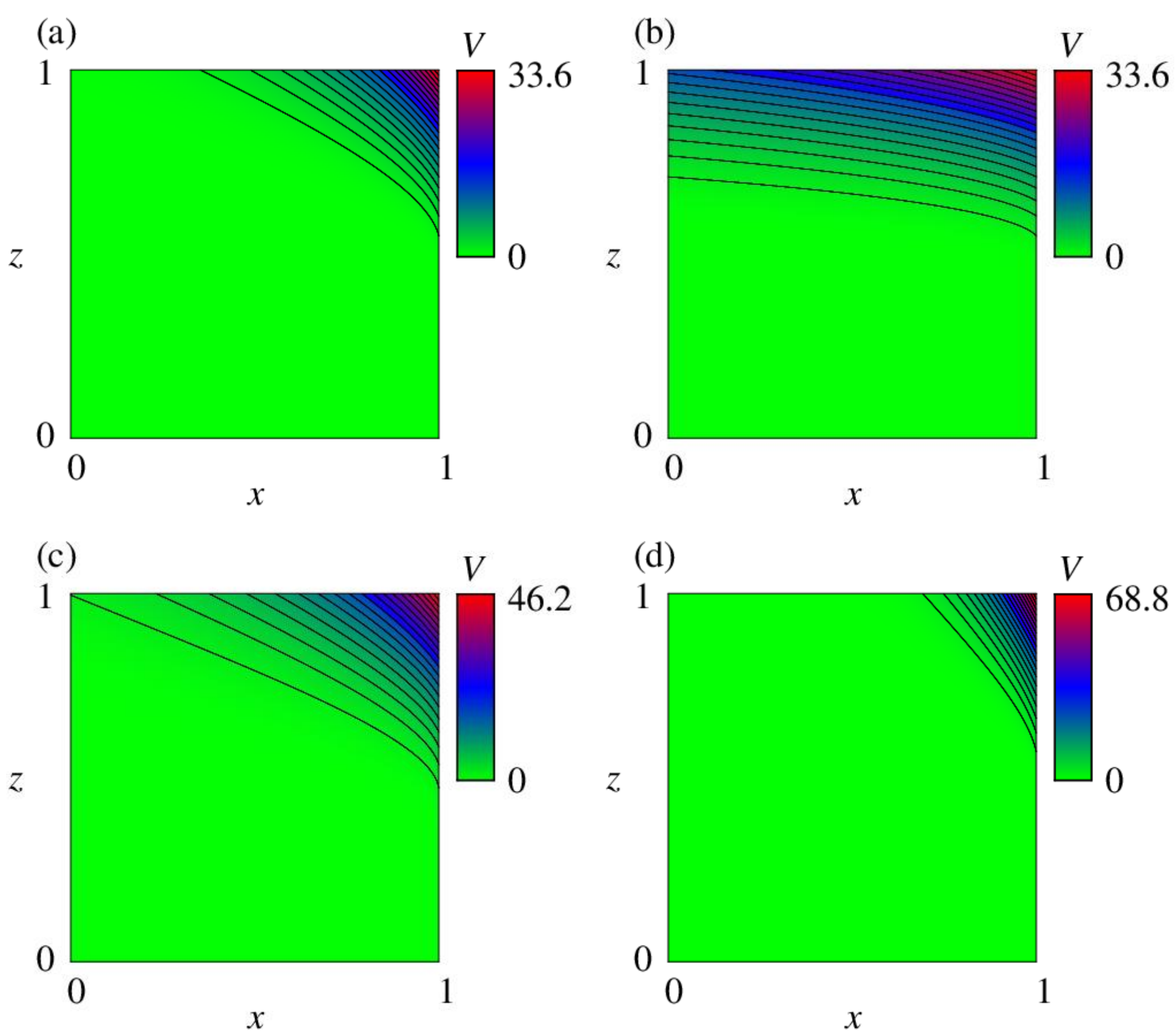

**Figure 4.** Computed $V$ with the filtered scheme for **(a)** the nominal case, **(b)** smaller $R = 0.05$, **(c)** larger $c = 0.6$, and **(d)** larger $\delta = 0.75$. The computational resolution is $N = 800$. Here, we set $f_2$.

### 4.3 Mean-field case

In this subsection, we analyze the nonlinear Kolmogorov equation. We investigate the two aspects concerning the present model and its discretization. The first is the discretization method (40) of the nonlinear term, and the second is the influence of the mean field effect. We set $\eta = 0$ and use $f_2$; hence, $V$ evaluates the regularized hitting probability unless otherwise specified. The computational resolution is $N = 800$. The computational analysis conducted here is exploratory; hence, a full theoretical analysis remains an open challenge.

Concerning the first issue, the computational results in the previous subsection suggest that $V$ is monotonically increasing with respect to both $x$ and $z$. This suggests that the same will apply to the nonlinear Kolmogorov equation, and we infer that the discretization method (40) simplifies to

$$\Lambda_{i,j} = 2\delta R\left(\frac{2\delta}{1-\delta}z_j + 1\right)z_j \times \frac{1}{2}\left(1+\tanh\left(\frac{V_{i,j}}{\kappa}\right)\right)\frac{V_{i,j}-V_{i,j-1}}{h}. \quad (51)$$

We set $\kappa = 0.5$ unless otherwise specified.

**Figure 5** shows the difference $\mathrm{Dify}_{i,j} = V_{i,j} - V_{i,j-1}$ for both monotone and fitted schemes, where we assume the second $f$. The computational results show that the obtained numerical solutions satisfy $\mathrm{Dify}_{i,j} \geq 0$ at each vertex, supporting their monotonically increasing property. Interestingly, the discretization with the filtered scheme satisfies this property even if it is not monotone by construction. We can therefore slightly simplify the discretization method (51) in practice. Although not related to the discretization method of the nonlinear term, we also examined $\mathrm{Difx}_{i,j} = V_{i,j} - V_{i-1,j}$, as shown in **Figure 6**, again yielding monotonicity.

We next investigate the parameter dependence of $V$ by focusing on $\kappa$ and $R$. We use the filtered scheme. **Figure 7** shows the parameter dependence with respect to $\kappa$ where $R = 0.20$ is fixed. **Figure 8** shows that for $R = 0.05$. For both values of $R$, increasing $\kappa$ leads to a less sharp profile of $V$, suggesting that the mean-field effect is smaller in these cases.

We also examine another function $\omega$ given by

$$\omega(v) = \frac{1}{2}\left(1+\tanh\left(\frac{v-1/2}{\kappa}\right)\right),\ v \in \mathbb{R}, \quad (52)$$

which has a sharp transition at $v = 1/2$. With this nonlinearity, the Kolmogorov equation is expected to have a sharp transition as well around region where $V = 1/2$, especially when $\kappa$ is small. A critical difference between the present and previous $\omega$ is that the former is bounded below by 0 but latter by 1/2; the present nonlinearity would therefore lead to a more drastic relative change in the solution to the Kolmogorov equation with the vanishing drift when $V$ is small. From a tourism management standpoint, this situation corresponds to a policy-maker who prefers adaptive as well as drastic changes in tourist flow to reduce overtourism.

**Figure 9** shows the numerical solutions to the nonlinear Kolmogorov equation with different values of $\kappa$ with (52). The numerical solution has a sharp transition, as expected for a small $\kappa$, creating an internal layer (see also **Figures 7-8** for comparison). Again, the numerical solutions do not have spurious oscillations, indicating that the filtered scheme works reasonably well.

We finally investigate the convergence of numerical solutions with $f = f_2, f_3$ to investigate their convergence speed depending on the regularity of the boundary data. Here, we regard the numerical solution with $N = 3200$ as the reference solution against which the convergence of numerical solutions with lower resolutions is discussed. We use $\eta = 0.1$ and examine the filtered scheme that has been suggested to be less diffusive in the previous subsection. **Tables 7-8** compare the $l^1$ (average absolute error at each vertex) and $l^\infty$ (maximum absolute error among all the vertices) errors between the numerical and reference solutions. For both $f_2$ and $f_3$, the convergence with respect to $l^1$ is slightly more optimistic than that with respect to $l^\infty$, but the difference is small; thus, in practice, using $l^\infty$ would suffice with the filtered scheme when one computes a similar Kolmogorov equation. We thus observed that the convergence is faster for the more regular $f$, as in the linear case.

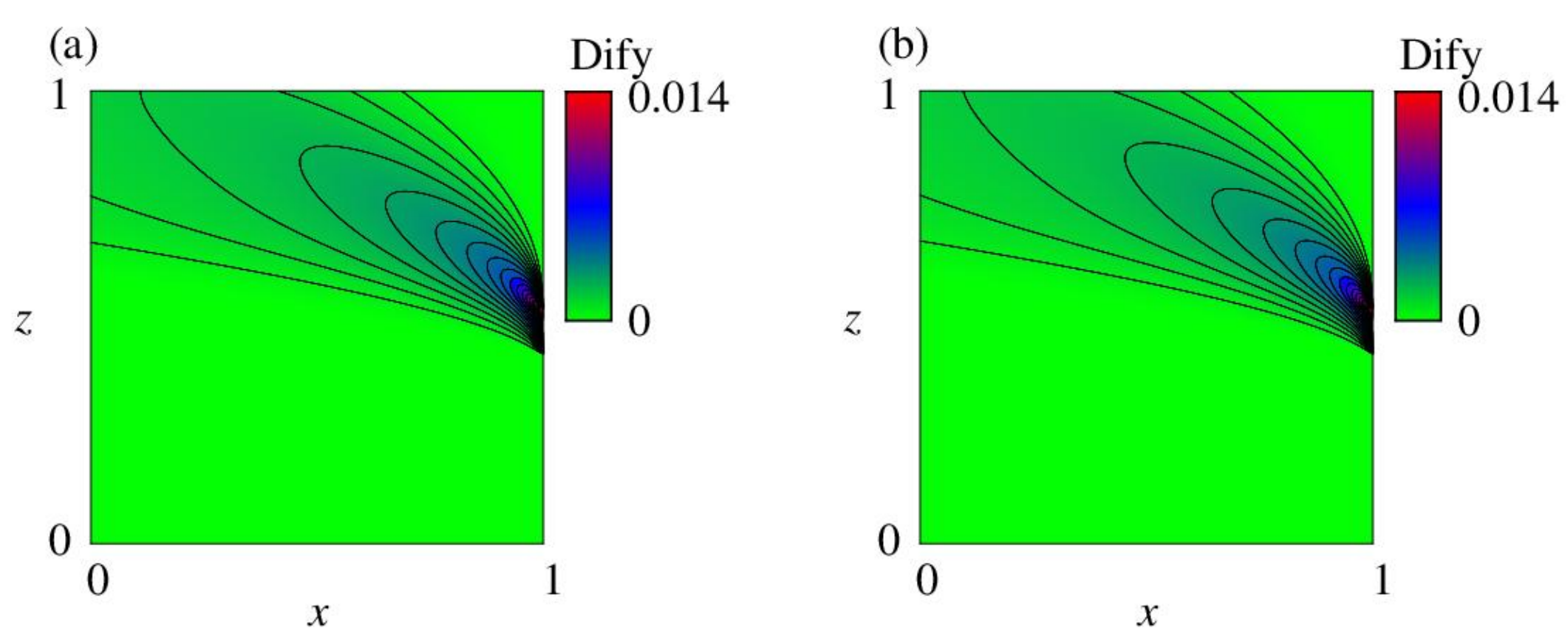

**Figure 5.** Difference $\mathrm{Dify}_{i,j} = V_{i,j} - V_{i,j-1}$ for **(a)** monotone and **(b)** filtered schemes.

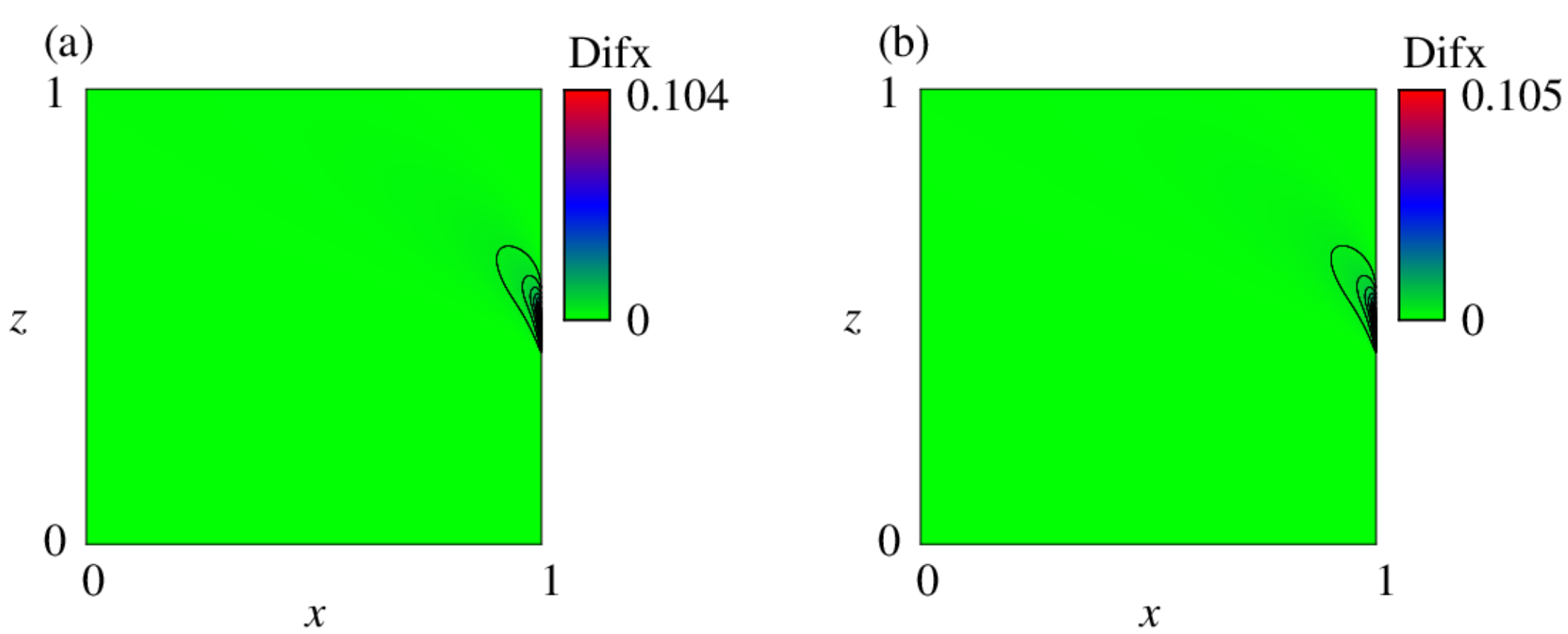

**Figure 6.** Difference $\mathrm{Difx}_{i,j} = V_{i,j} - V_{i-1,j}$ for **(a)** monotone and **(b)** filtered schemes.

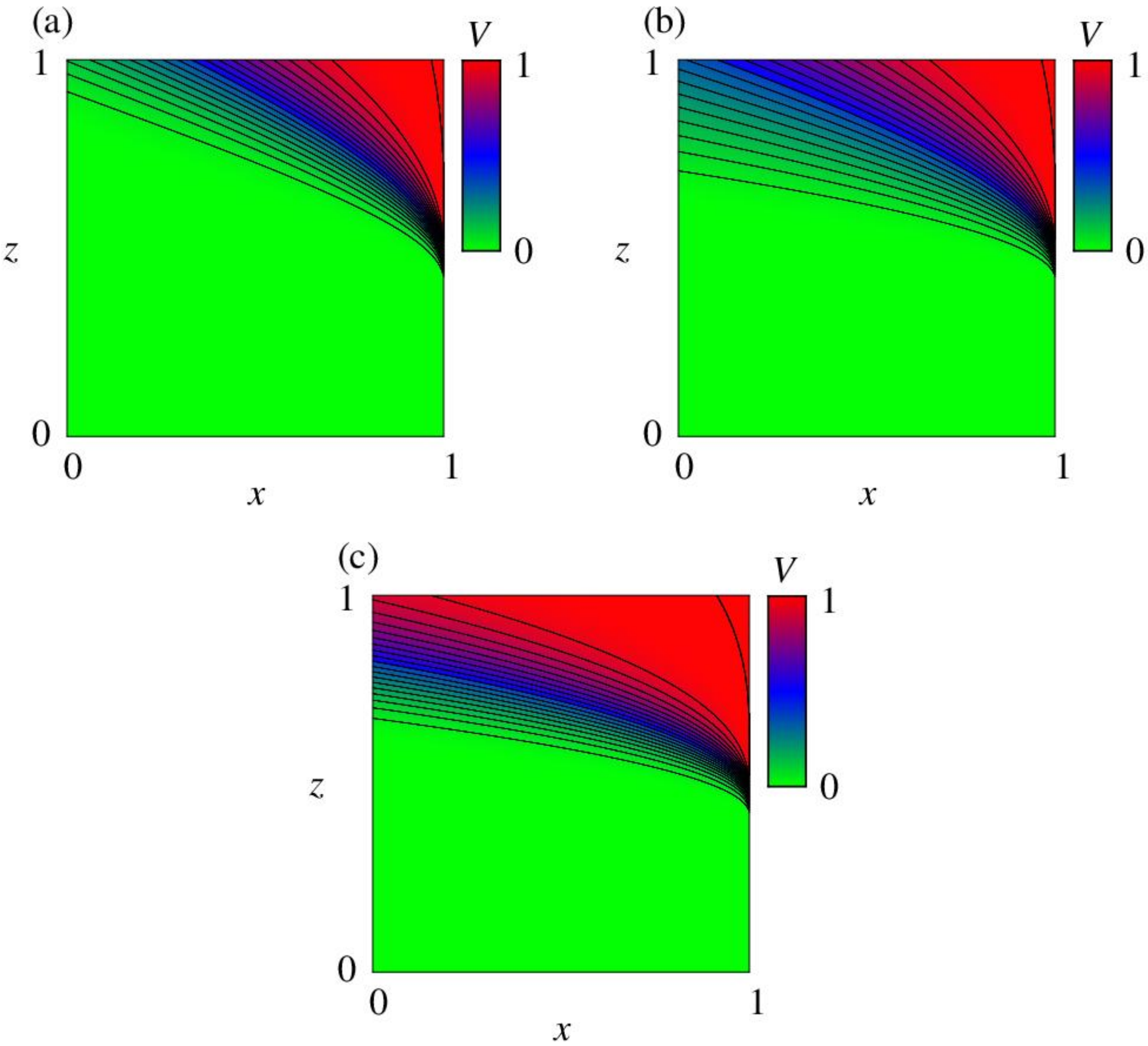

**Figure 7.** Parameter dependence with respect to $\kappa$: **(a)** $\kappa = 0.0001$, **(b)** $\kappa = 0.5$, and **(c)** $\kappa = 5$. Here, we set $R = 0.20$.

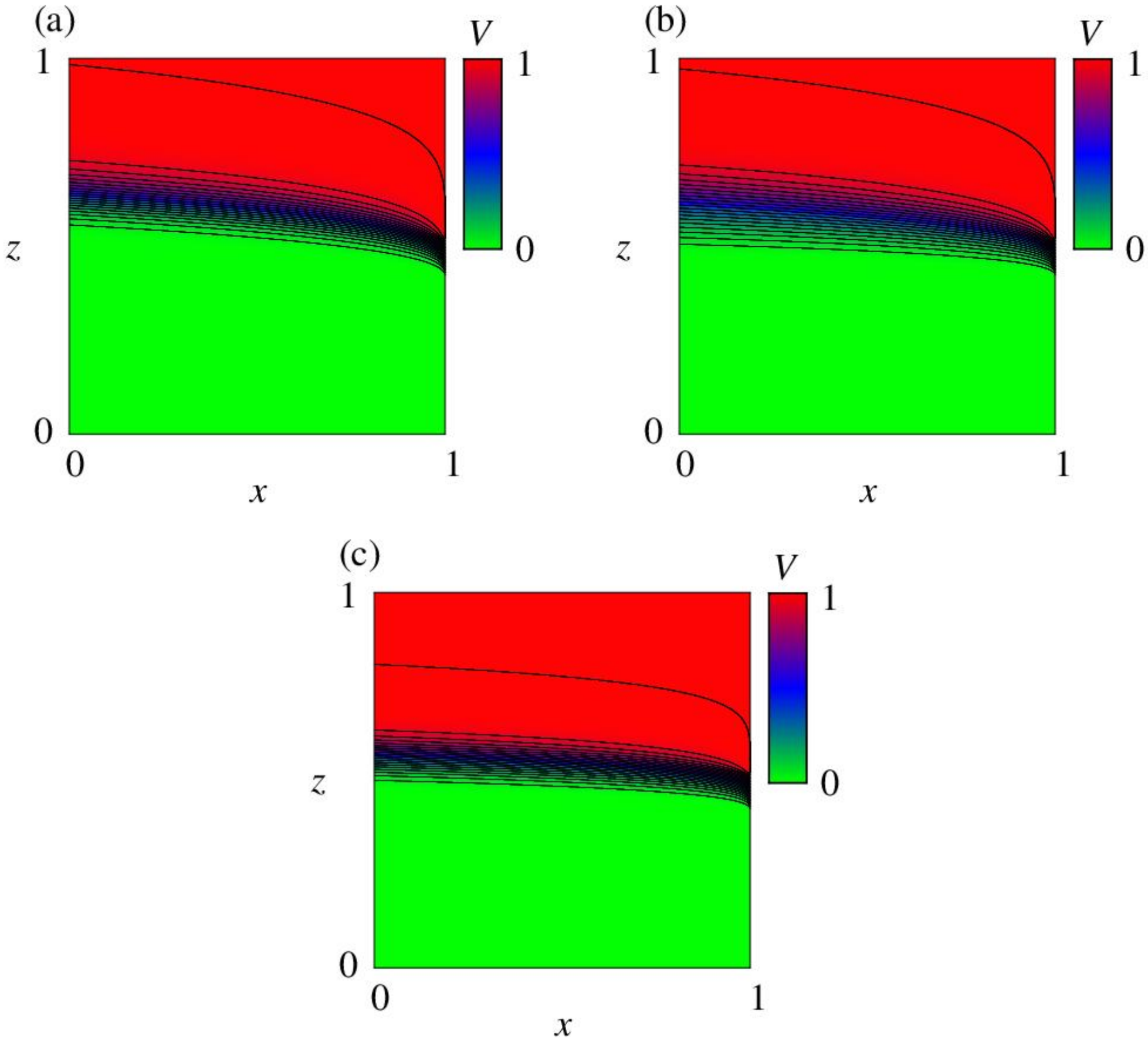

**Figure 8.** Parameter dependence with respect to $\kappa$: **(a)** $\kappa = 0.05$, **(b)** $\kappa = 0.5$, and **(c)** $\kappa = 5$. Here, we set $R = 0.05$.

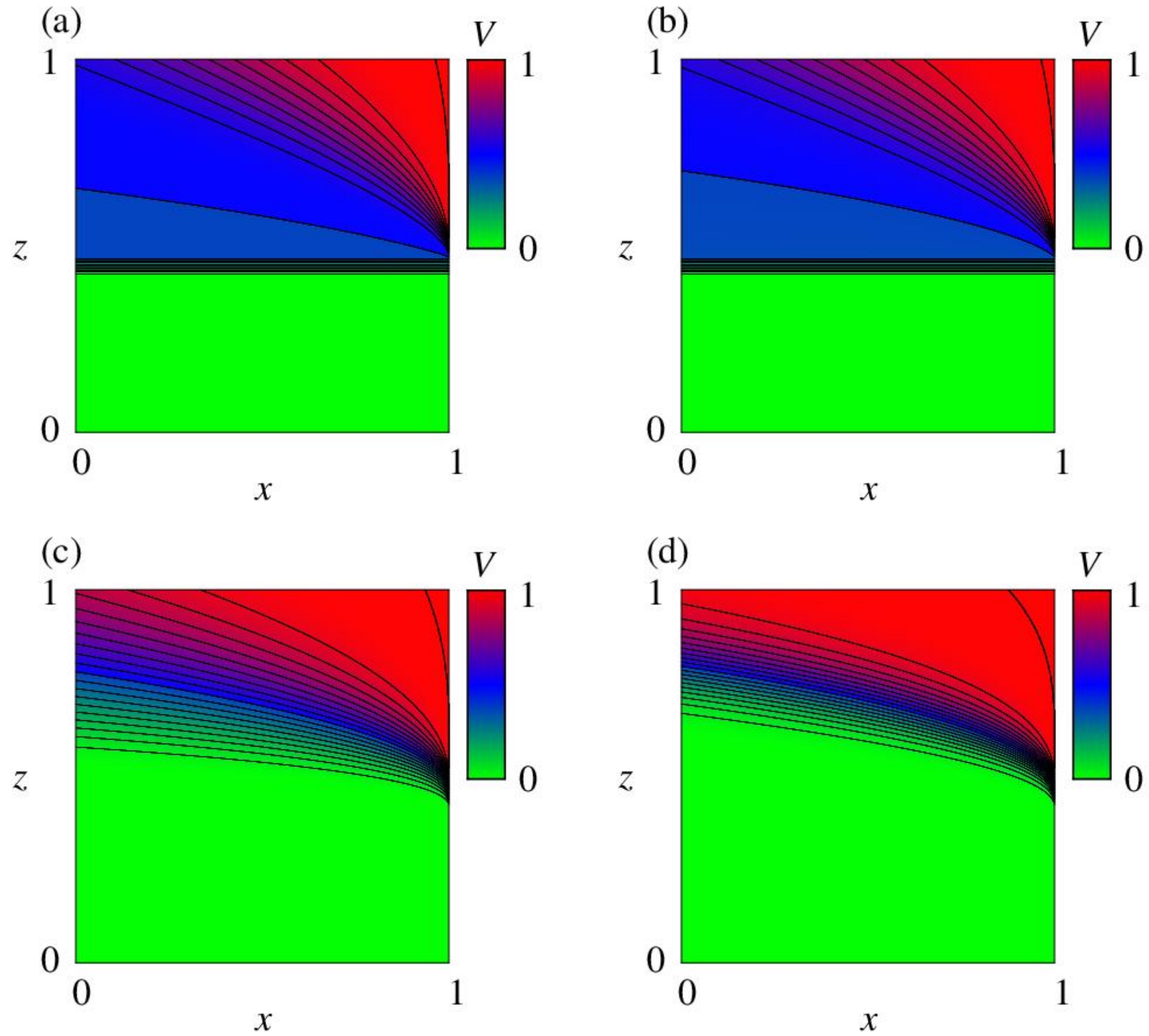

**Figure 9.** Parameter dependence with respect to $\kappa$ for the nonlinearity (52): **(a)** $\kappa = 0.0001$, **(b)** $\kappa = 0.01$, **(c)** $\kappa = 1$, and **(c)** $\kappa = 100$. Here, we set $R = 0.2$.

**Table 7.** Comparison between the numerical and reference solutions: $f = f_2$.

| $N$ | $l^1$ | | $l^\infty$ | |
|---|---|---|---|---|
| | Error | Convergence rate | Error | Convergence rate |
| 200 | 0.00333 | 0.783 | 0.06413 | 0.563 |
| 400 | 0.00194 | 0.864 | 0.04342 | 0.728 |
| 800 | 0.00105 | 1.256 | 0.02622 | 1.121 |
| 1600 | 0.00045 | | 0.01206 | |

**Table 8.** Comparison between the finite difference and reference solutions: $f = f_3$.

| $N$ | $l^1$ | | $l^\infty$ | |
|---|---|---|---|---|
| | Error | Convergence rate | Error | Convergence rate |
| 200 | 0.04179 | 1.131 | 0.39855 | 1.055 |
| 400 | 0.01909 | 1.300 | 0.19181 | 1.239 |
| 800 | 0.00775 | 1.576 | 0.08125 | 1.486 |
| 1600 | 0.00260 | | 0.02900 | |

## 5. Conclusion

We proposed Jacobi diffusion driven by unsteady drift and analyzed its well-posedness. Its mean-field version was also introduced. The Kolmogorov equations associated with these processes were derived, and their properties, including stochastic representations, were discussed. These equations require boundary conditions only along a part of the boundary because of the directions of characteristics and degenerate diffusion. A finite difference method equipped with a filtered scheme was presented to discretize both linear and nonlinear Kolmogorov equations in a unified manner. This finite difference method generates a unique solution for both linear and nonlinear Kolmogorov equations. The computational results revealed that the regularity of boundary conditions critically affects the convergence rate of the finite difference method. Moreover, the proposed finite difference method with the filtered scheme proved effective for the Kolmogorov equations, as it generated numerical solutions without spurious oscillations.

A major limitation of this paper is that the results regarding the uniqueness of the solution remain limited, in particular. Although the existence of a solution to the first Kolmogorov equation was shown, an unresolved issue is the assumption made in the comparison theorem, which is conditional unless the ordering property of subsolution and supersolution is proven. A critical assumption in this paper is that the driving process is deterministic, although the system driven by a stochastic process is more general, which is an important direction for future research. A potential issue in this case is the loss of computational simplicity in numerically solving the Kolmogorov equations because our finite difference method exploits the absence of diffusion in the governing equation of $Y$. Moreover, the boundary condition will have to be reconsidered, especially if both boundaries $x=0,1$ become touchable by $X$. Nevertheless, the filtered scheme can also be applied to such cases and would be beneficial if the boundary condition is sufficiently smooth. We are currently investigating diffusion processes, including the Jacobi process, whose coefficients are modulated by long-term memory effects, resulting in a non-Markovian process. Such models can be particularly relevant in water environmental management problems. Another future direction for this study includes investigating nonlinear Kolmogorov equations where the drift and/or diffusion coefficients in the $x$ direction depend on $V$. This case is more challenging because the Fichera function depends on $V$ itself.

## Appendices

### A.1 Proofs

***Proof of Proposition 1***

The dynamics of $Z$ are decoupled from those of $X$; hence, there exists a unique strong solution to the second equation of the system (9) (i.e., a logistic curve). With respect to the process $X$, to apply Theorem 1 in Yamada and Watanabe [55] (YW1971 in the sequel), we introduce the auxiliary system in which the coefficient $\sqrt{X_t(1-X_t)}$ in the diffusion term is replaced by $\sqrt{|X_t(1-X_t)|}$ and $X_t$ and $Z_t$ in the drift term by $\max\{0,\min\{X_t,1\}\}$ and $\max\{0,\min\{Z_t,1\}\}$, respectively:

$$\mathrm{d}\begin{pmatrix} X_t \\ Z_t \end{pmatrix} = \begin{pmatrix} 2\delta\max\{0,\min\{Z_t,1\}\}+1-\delta-\max\{0,\min\{X_t,1\}\} \\ -2\delta R\left(\dfrac{2\delta}{1-\delta}\max\{0,\min\{Z_t,1\}\}+1\right)\max\{0,\min\{Z_t,1\}\} \end{pmatrix}\mathrm{d}t + \begin{pmatrix} c\sqrt{|X_t(1-X_t)|} \\ 0 \end{pmatrix}\mathrm{d}B_t,\ t>0. \quad (53)$$

The auxiliary system admits at most one strong solution because the assumption made in Theorem 1 of YW1971 is satisfied because of the Lipschitz continuity of the drift coefficient, the Hölder continuity of the diffusion coefficient, and the diagonal diffusion coefficient (we can take $\rho(x)=Cx^{1/2}$ with a constant $C>0$ for the function $\rho$ used in this literature). Here, the existence up to time $\tau$ is due to the boundedness and local Lipschitz continuity of the drift and diffusion terms along with the suitable regularization as employed in Proof of theorem 6.1.1 in Alfonsi [18] (see also Theorem 6.5 in Chapter 2 of Mao [65]). Then, the strong solution to the auxiliary system (53) coincides with that of the original one (9) for the time interval $(0,\tau)$. Owing to the local Lipschitz continuity of the drift and diffusion terms of the original system (9), this solution is unique to this system up to time $\tau$. Note that we already know $Z_t=\max\{0,\min\{Z_t,1\}\}$ for $t\geq 0$, and hence, taking $\max\{0,\min\{\cdot,1\}\}$ is innocuous. Moreover, we have $\max\{0,\min\{X_t,1\}\}=X_t$ and $\sqrt{|X_t(1-X_t)|}=\sqrt{X_t(1-X_t)}$ for $0<t<\tau$.

□

***Proof of Proposition 2***

We can use the comparison result of the SDEs by viewing the system (53) as one-dimensional, where $Z$ is given by a smooth and bounded curve. We use Theorem 1.1 of Ikeda and Watanabe [56] (IW1977 in the sequel); we take $\rho(x)=Cx^{1/2}$ with a constant $C>0$ for the function $\rho$ used in this literature. We consider a process $W=(W_t)_{t\geq 0}$ governed by

$$\mathrm{d}W_t=(1-\delta-W_t)\mathrm{d}t+c\sqrt{W_t(1-W_t)}\mathrm{d}B_t,\ t>0 \quad (54)$$

subject to the initial condition $W_0=X_0$. The process $W$ is a classical Jacobi diffusion and is bounded in $D$ because $\delta\in(0,1)$ and $2(1-\delta)>c^2$ (the left most inequality of (7)). We can apply Theorem 1.1 of

IW1977 to (53) and (54), showing that $W_t \le X_t$ for all $t \in [0,\tau)$. Indeed, we have

$$\underbrace{2\delta \max\{0,\min\{z,1\}\} + 1 - \delta - x}_{\text{Drift of the original system}} - \underbrace{(1-\delta-x)}_{\text{Drift of the classical Jacobi diffusion}} = 2\delta \max\{0,\min\{z,1\}\} > 0 \tag{55}$$

for $x \in \mathbb{R}$ and $z > 0$. Owing to $W_t > 0$ a.s. $t > 0$ (Chapter 6 in Alfonsi [18]), we obtain $0 < W_t \le X_t$ a.s. $t > 0$, and hence, the desired result (10). This, combined with **Proposition 1**, completes the proof.

□

***Proof of Proposition 3***

**Proposition 2** shows that the hitting of the process $X$ to the lower boundary never occurs. We prepare an upper-bounding stochastic process from above so that, again, Theorem 1.1 of IW1977 [56] can be used. The strategy of the proof is therefore qualitatively the same as that of **Proposition 2**, but now we use a different auxiliary process $W$. We consider the process $W = (W_t)_{t\ge0}$ as the unique solution to the SDE

$$\mathrm{d}W_t = (2\delta\rho + 1 - \delta - W_t)\mathrm{d}t + c\sqrt{W_t(1-W_t)}\mathrm{d}B_t,\ t > 0 \tag{56}$$

subject to the initial condition $W_0 = X_0$. This $W$ is a classical Jacobi process. We have

$$2\delta\rho + 1 - \delta - W_t = 2\delta\left(\frac{1}{2} - \frac{c^2}{4\delta}\right) + 1 - \delta - W_t = 1 - \frac{c^2}{2} - W_t,\ t > 0. \tag{57}$$

We have $Z_0 \le \rho$, and because $Z_t$ is monotonically decreasing from 1 to 0 with respect to $t > 0$, it follows that $Z_t \le \rho$ for $t \ge t_c$. Therefore, if $z < \rho$, then we have

$$\underbrace{2\delta\rho + 1 - \delta - W}_{\text{Drift of the classical Jacobi diffusion}} - \underbrace{(2\delta \max\{0,\min\{z,1\}\} + 1 - \delta - x)}_{\text{Drift of the original system}} = 2\delta(\rho - \max\{0,\min\{z,1\}\}) > 0. \tag{58}$$

Because $W_t < 1$ a.s. $t > 0$ (Chapter 6 in Alfonsi [18]), we obtain $X_t \le W_t < 1$ a.s. $t > 0$ under the assumption $Z_0 \le \rho$, and hence, the desired result (10).

□

***Poof of Proposition 4***

In this proof, we abbreviate Rokhlin [37] as RO14. A difference between RO14 and this paper is that the former assumes the Lipschitz continuity of drift and diffusion coefficients in a domain, whereas the latter is a non-Lipschitz diffusion coefficient. Lipschitz continuity was assumed in RO14 to justify the pathwise uniqueness of the SDE. We proved the pathwise continuity of our target system in **Proposition 2**, overcoming this difficulty. Another difference is that RO14 assumes a control problem, whereas ours does not, implying that our method is simpler in this sense. By using notations in RO14, we can obtain the statement by setting $G = D \times D$, $\bar{G} = \bar{D} \times \bar{D}$, and $\hat{G} = (\bar{D} \times \bar{D}) \setminus \{z = 0,1\}$, with which a suitable stochastic subsolution $\underline{V}$ and supersolution $\bar{V}$ can be constructed according to Theorems 2 (with a null control) and Theorem 3 of RO14. The fact that $V$ is a viscosity solution to the Kolmogorov equation (13) follows from the dynamic programming principle by using $\underline{V}$ and $\overline{V}$ as in Remark 4 of RO14.

□

***Proof of Proposition 5***

This is a direct application of Theorem 1 of Rokhlin [37] to our case along with **Proposition 4**.

□

***Proof of Proposition 6***

According to Theorem 8 of Oberman [51], it suffices to show the following two inequalities at each vertex exterior to $\Gamma_{\mathrm{num}}$ (those at each vertex in $\Gamma_{\mathrm{num}}$ is trivial):

$$\frac{\partial \Lambda_{i,j}}{\partial V_{i,j}}, -\frac{\partial \Lambda_{i,j}}{\partial V_{i,j-1}} \geq 0 \tag{59}$$

for $i, j = 0,1,2,...,N$. Note that $\Lambda_{i,j}$ is Lipschitz continuous with respect to both $V_{i,j}$ and $V_{i,j-1}$, which can be directly checked from the discussion below along with the elementary results for any $v \in \mathbb{R}$:

$$\left|\tanh\left(\frac{v}{\kappa}\right)\right| \leq 1 \text{ and } \left|\cosh^2\left(\frac{v}{\kappa}\right)\right| \leq C_\kappa v \text{ with a constant } C_\kappa > 0. \tag{60}$$

The first inequality of (59) is proven as follows: if $V_{i,j} - V_{i,j-1} \geq 0$, then

$$\begin{aligned} \frac{\partial \Lambda_{i,j}}{\partial V_{i,j}} &= 2\delta R\left(\frac{2\delta}{1-\delta} z_j + 1\right) z_j \times \frac{\partial}{\partial V_{i,j}}\left\{\frac{1}{2}\left(1+\tanh\left(\frac{V_{i,j}}{\kappa}\right)\right)\frac{V_{i,j} - V_{i,j-1}}{h}\right\} \\ &= 2\delta R\left(\frac{2\delta}{1-\delta} z_j + 1\right) z_j \left\{\frac{1}{2\kappa}\frac{V_{i,j} - V_{i,j-1}}{h}\left[\cosh\left(\frac{V_{i,j}}{\kappa}\right)\right]^{-2} + \frac{1}{2}\left(1+\tanh\left(\frac{V_{i,j}}{\kappa}\right)\right)\frac{1}{h}\right\} \\ &\geq 0 \end{aligned} \tag{61}$$

and if $V_{i,j} - V_{i,j-1} < 0$ then

$$\begin{aligned} \frac{\partial \Lambda_{i,j}}{\partial V_{i,j}} &= 2\delta R\left(\frac{2\delta}{1-\delta} z_j + 1\right) z_j \times \frac{\partial}{\partial V_{i,j}}\left\{\frac{1}{2}\left(1+\tanh\left(\frac{V_{i,j-1}}{\kappa}\right)\right)\frac{V_{i,j} - V_{i,j-1}}{h}\right\} \\ &= 2\delta R\left(\frac{2\delta}{1-\delta} z_j + 1\right) z_j \frac{1}{2}\left(1+\tanh\left(\frac{V_{i,j-1}}{\kappa}\right)\right)\frac{1}{h} \\ &\geq 0 \end{aligned} . \tag{62}$$

The second inequality of (59) is proven as follows: if $V_{i,j} - V_{i,j-1} \geq 0$, then

$$\begin{aligned} \frac{\partial \Lambda_{i,j}}{\partial V_{i,j-1}} &= 2\delta R\left(\frac{2\delta}{1-\delta} z_j + 1\right) z_j \times \frac{\partial}{\partial V_{i,j-1}}\left\{\frac{1}{2}\left(1+\tanh\left(\frac{V_{i,j}}{\kappa}\right)\right)\frac{V_{i,j} - V_{i,j-1}}{h}\right\} \\ &= -2\delta R\left(\frac{2\delta}{1-\delta} z_j + 1\right) z_j \frac{1}{2}\left(1+\tanh\left(\frac{V_{i,j}}{\kappa}\right)\right)\frac{1}{h} \\ &\leq 0 \end{aligned} \tag{63}$$

and if $V_{i,j} - V_{i,j-1} < 0$ then

$$\frac{\partial \Lambda_{i,j}}{\partial V_{i,j-1}} = 2\delta R\left(\frac{2\delta}{1-\delta}z_j+1\right)z_j \times \frac{\partial}{\partial V_{i,j-1}}\left\{\frac{1}{2}\left(1+\tanh\left(\frac{V_{i,j-1}}{\kappa}\right)\right)\frac{V_{i,j}-V_{i,j-1}}{h}\right\}$$

$$= 2\delta R\left(\frac{2\delta}{1-\delta}z_j+1\right)z_j\left\{\frac{1}{2\kappa}\left[\cosh\left(\frac{V_{i,j-1}}{\kappa}\right)\right]^{-2}\frac{V_{i,j}-V_{i,j-1}}{h}-\frac{1}{2}\left(1+\tanh\left(\frac{V_{i,j-1}}{\kappa}\right)\right)\frac{1}{h}\right\}. \quad (64)$$

$$\leq 0$$

The proposition is proven because (61)-(64) lead to (59).

□

### A.2 On the continuity assumption of Proposition 5

The main conditions assumed in Theorem 5.1 of Koike [60] to apply a comparison argument to our Kolmogorov equation are equations (5.3), (5.4), and (5.6) in this literature. The condition (5.3) is satisfied in our case because we are dealing with a linear Kolmogorov equation whose drift and diffusion coefficients are Lipschitz continuous on $\bar{D}\times\bar{D}$. Finally, condition (5.6) is satisfied in our case because our domain is a rectangle.

What is problematic in our case is condition (5.4), which in our context is stated as follows: There exists a nonnegative, continuous, and nondecreasing function $\varsigma:[0,+\infty)\to[0,+\infty)$ with $\varsigma(0)=0$ such that if $M_1, M_2 \in \mathbb{R}$ and $\mu>1$ satisfy the inequality

$$-3\mu\left(w_1^2+w_2^2\right)\leq w_1^2 X_1 - w_2^2 X_2 \leq 3\mu\left(w_1-w_2\right)^2 \text{ for any } w_1, w_2\in\mathbb{R} \quad (65)$$

then

$$F\left(x_2,z_2,p,q,M_2\right)-F\left(x_1,z_1,p,q,M_1\right)\leq\varsigma\left(\left\|\mathbf{x}_1-\mathbf{x}_2\right\|\left(1+|p|+|q|+\mu\left\|\mathbf{x}_1-\mathbf{x}_2\right\|\right)\right) \quad (66)$$

with $\left\|\mathbf{x}_1-\mathbf{x}_2\right\|=|x_1-x_2|+|z_1-z_2|$ for any $x_1,x_2,z_1,z_2\in\bar{D}$ and $p_1,p_2\in\mathbb{R}$, where

$$F\left(x,z,p,q,M\right)=-\left(2\delta z+1-\delta-x\right)p+2\delta R\left(\frac{2\delta}{1-\delta}z+1\right)zq-\frac{1}{2}c^2x\left(1-x\right)M . \quad (67)$$

Now, we check the statement above to show that it fails in our case. For each $x_1,x_2,z_1,z_2\in\bar{D}$, $M_1,M_2\in\mathbb{R}$, and $\mu>1$, we set

$$w_1=\sqrt{\frac{c^2}{2}x_1\left(1-x_1\right)} \text{ and } w_2=\sqrt{\frac{c^2}{2}x_2\left(1-x_2\right)} \text{ for any } w_1,w_2\in\mathbb{R} \quad (68)$$

to specify (65) as follows:

$$\begin{aligned}-\frac{3\mu c^2}{2}\left\{x_1\left(1-x_1\right)+x_2\left(1-x_2\right)\right\}&\leq\frac{c^2}{2}x_1\left(1-x_1\right)M_1-\frac{c^2}{2}x_2\left(1-x_2\right)M_2\\&\leq\frac{3\mu c^2}{2}\left(\sqrt{x_1\left(1-x_1\right)}-\sqrt{x_2\left(1-x_2\right)}\right)^2\end{aligned}. \quad (69)$$

Then, by (67) and (69), the left-hand side of (66) is rewritten as follows:

$$F\left(x_2,z_2,p,q,M_2\right)-F\left(x_1,z_1,p,q,M_1\right)\leq\bar{C}\left\|\mathbf{x}_1-\mathbf{x}_2\right\|\left(|p|+|q|\right)+\frac{3\mu c^2}{2}\left(\sqrt{x_1\left(1-x_1\right)}-\sqrt{x_2\left(1-x_2\right)}\right)^2. \quad (70)$$

Here, $\bar{C} > 0$ is a constant depending only on $\delta$ and $R$. Owing to (70), it is impossible to find a suitable right-hand side of (66) to bound the last line of (70) from above. This difficulty comes from the non-Lipschitz continuity of the diffusion coefficient to which the existing comparison principles and continuity results do not apply [66,67].

### A.3 Computational cost

We study the average iteration count to obtain numerical solutions among all $j$ (**Table A1**) for each computational case with $f = f_2$. The average iteration count increases as $N$ increases, and interestingly, it is greater in the monotone scheme than in the filtered scheme; their difference is approximately 20%. The only difference between the two schemes is the discretization of the first-order partial differential terms; hence, the results obtained suggest that the use of the filtered scheme performs better, as shown in **Figure 3**. Finally, the best least-squares proportional estimates of the average interaction counts for the monotone and filtered schemes as functions of $N$ are $1.572N$ ($R^2$=0.99) and $1.274N$ ($R^2$=0.99), respectively.

**Table A1.** Average iteration count for the monotone and filtered schemes.

| $N$ | Monotone | Filtered |
|---|---|---|
| 200 | 455.5 | 355.6 |
| 400 | 785.2 | 629.9 |
| 800 | 1412.3 | 1177.1 |
| 1600 | 2613.1 | 2327.4 |
| 3200 | 4913.4 | 3873.1 |